\documentclass[12pt,thmsa,twoside,a4paper]{article}
\pdfoutput=1

\usepackage{graphicx}
\usepackage{ifthen,calc}
\usepackage{epstopdf}
\usepackage{rotating}
\usepackage{pdflscape}
\usepackage{longtable}
\usepackage{subfig}
\usepackage{float}
\usepackage{epstopdf}
\usepackage{xcolor}

\newcommand\bovermat[2]{%
  \makebox[0pt][l]{$\smash{\overbrace{\phantom{%
    \begin{matrix}#2\end{matrix}}}^{\text{#1}}}$}#2}

\usepackage{multicol} 
\usepackage[reqno]{amsmath}
\usepackage{amssymb} 
\usepackage{amscd} 
\usepackage{amsfonts}

\usepackage{bm}       
\usepackage{makeidx}
\usepackage{amsxtra}  
\usepackage{epsf}     
\usepackage{pifont}

\newcommand{\PreserveBackslash}[1]{\let\temp=\\#1\let\\=\temp}

\newcommand{\ThesisBibNote}[1]{}

\newcommand{\eg}{e.g.,  }

\newcommand{\IndexEmphType}{\textbf{Bold}}       

\newcommand{\I}[1]{#1\index{#1|IndexEmph}}

\makeatletter 

                 { \clearpage
                   \setlength{\parindent}{0mm}
                   \addcontentsline{toc}{part}{Index}%
                   \pagestyle{myheadings}
                   \markboth{Index}{Index}
                     \let\item\@idxitem
                   \begin{multicols}%
                         {2}[{%
                              \section*{\indexname}
                              The list of \emph{symbols} is sorted
                              alphabetically by the names of the
                              defined symbols.
                              For example, all kinds of brackets are
                              sorted under ``b'', $\Hermitian{A}$
                              under ``s'' (as in star), $\overline{a}$
                              under ``o'' (as in overline) and $\dual{g}$ under
                              ``t'' (as in tilde).
                              \IndexEmphType\ page numbers are used to
                              indicate pages with important information
                              about the entry, \eg a precise definition
                              or the most detailed explanation found in this
                              thesis, while page
                              numbers in normal type indicate a textual
                              reference.
                             }
                            ]
                   \par\bigskip%
                 }%
                 {\end{multicols}}

\makeatother  

\usepackage[swedish,english]{babel}
\renewcommand{\vec}[1]{{\bm{#1}}} 

\newcommand{\dual}[1]{\widetilde{#1}}

\newcommand{\Hermitian}[1]{#1^*}

\newcommand{\mat}[1]{#1}  

\newcommand{\vecphi}{\bm{\phi}}

\newcommand{\defequal}{\ensuremath{\stackrel{\mathrm{def}}{=}}}

\newcommand{\pdif}[2]{\frac{\partial #1}{\partial #2}}

\usepackage[pdftex,hyperindex]{hyperref}  

   \newcommand{\e}{\protect\ensuremath{\mathrm{e}}} 

\usepackage{amsthm} 

\theoremstyle{plain}

\theoremstyle{definition}

\theoremstyle{remark}

\usepackage{moreverb}

  \cleardoublepage

\makeindex
\usepackage{vmargin}
\setpapersize{A4}
\setmarginsrb{2cm}{8mm}{2cm}{3cm}{0mm}{12mm}{0mm}{0mm}

\usepackage{varioref}

\usepackage{afterpage}
\allowdisplaybreaks

\begin{document}

\author{Niklas Grip$^1$%
        \thanks{\scriptsize The authors were supported by
                the \href{http://www.formas.se/default____529.aspx}{Swedish Research Council Formas} grants
                (registration numbers 2007--1430 and 2012--1037) as well as by the Swedish Construction Industry's Organisation for Research and Development (SBUF) grant 13010.}
        , 
        Natalia Sabourova$^{1*,}$
        \addtocounter{footnote}{5}%
        \thanks{\scriptsize The author was supported by Elsa and Sven Thysells Foundation for Structural Engineering Studies at Lule{\aa} University of Technology.}
        ~
        and Yongming Tu$^2$
        \thanks{\scriptsize The author was supported by the National Natural Science Foundation of China (project number 51378104).}
        \\
        {\scriptsize
         $^1$
         Lule{\aa} University of Technology,
         SE-971 87 Lule{\aa},
         Sweden,
         \href{mailto:Niklas.Grip@ltu.se}{Niklas.Grip@ltu.se},
         \href{mailto:Natalia.Sabourova@ltu.se}{Natalia.Sabourova@ltu.se}.
        }
        \\
        {\scriptsize
         $^2$
         School of Civil Engineering, Southeast University, Nanjing, China,
         \href{mailto:tuyongming@seu.edu.cn}{tuyongming@seu.edu.cn},\href{mailto:yongming.tu@ltu.se}{yongming.tu@ltu.se}.
        }%
        }
\title{
       Sensitivity-Based Model Updating  \\
       for Structural Damage Identification  \\
       Using Total Variation Regularization
       }
\date{\vspace{-5ex}}
\maketitle
\thispagestyle{empty}
\begin{abstract}
Sensitivity-based Finite Element Model Updating (FEMU) is one of the widely accepted techniques used for damage identification
in structures. FEMU can be formulated as a numerical optimization problem and solved iteratively making automatic updating of
the uncertain model parameters by minimizing the difference
between measured and analytical structural properties. However, in the presence of noise in the measurements,
the updating results are usually prone to errors. This is mathematically described as instability of the damage identification as an inverse problem. One way to resolve this problem is by using regularization.
In this paper we investigate regularization methods based on the minimization of the total variation of the uncertain
model parameters and compare this solution with a rather frequently used regularization based on an interpolation technique. For well-localized damages the results show a clear
advantage of the proposed solution in terms of the identified location and severity of damage compared with the interpolation based solution.

For a practical test of the proposed method we use a reinforced concrete plate. Measurements and analysis were
repeated first on an undamaged plate, and then after applying four different
degrees of damage.
\\
%
%
  %
  \textbf{Keywords:}
  Finite element model updating, damage identification, total variation regularization, (pseudo) Huber function, interpolation, reinforced concrete plate
\end{abstract}


\section{Introduction}
In this paper, we deal with finite element model updating by the classical iterative sensitivity based method \cite{DFP98SHM,MLF11SHM}. Compared to other finite element model updating methods, the sensitivity based method showed computational efficiency and good sensitivity to small damages
\cite{Mar10SHM,MLF11SHM}. Basically, there are two application areas of model updating. In the first place, it is applied in order to increase the reliability of the finite element model and thus, for example, the prediction of the dynamic behavior of the structure under different loads. Another application area is damage identification in structures which is the focus of this paper.


In mathematical language, damage identification by finite element model updating  is a parameter estimation problem. The finite element model is parameterized by uncertain parameters, which are updated by some parameter estimation technique. We assume here that the model is physically meaningful and thus can accurately represent the behavior of the actual structure, so that the damage identification problem can be reduced to the parameter estimation only.
The parameter estimation problems belong to a class of inverse problems, i.e. knowing the model outputs, one need to obtain the internal model parameters. In the presence of noise in the outputs, which is the case with vibration tests, the inverse problem becomes ill-posed, i.e. small variations in the outputs lead to unreasonably large variations in the model parameters. Such problems can be solved by using regularization, which is increasingly more often consistently taken into account in the area of structural damage identification (\cite{Natke88SHM,FrMo95SHM,Link01SHM,FMA01SHM,TMR02SHM,Weber09SHM}, etc).

In this paper, we investigate a regularization tool for the ill-posed damage identification problem that has its origin in
image processing and which is associated with minimization of the total variation of the uncertain parameters.
We compare then this regularization technique with  a rather frequently used  interpolation
with so-called damage functions introduced in \cite{TMR02SHM}.
By using damage or interpolating functions, the algorithm is free
to choose parameter values freely on a more sparse grid,
and then parameters in the intermediate points are chosen by interpolation. This gives smoothing, but you lose resolution.
It would be desirable to keep the highest possible resolution but add restrictions to the damage identification algorithm
that favors solutions with a sharp increase of damage index close to a damage and keep damage index close to zero elsewhere.
We show that the total variation regularization brings the parameter estimation close to the desirable solution and in the case
of well-localized damage, it results in a more precise damage identification than the interpolation method.

In reinforced concrete structures shear cracks may form well-localized damage patterns. When such cracks develop, a brittle failure of the structure may be close --- an inclined crack can find its way through a structure, without being prevented by reinforcement. It is of great interest to identify location and severity of such local cracks more precisely without smoothing the damage to the areas nearby and in this way to distinguish these cracks from other less severe cracks, such as e.g. bending cracks.

\subsection{Damage parametrization}
A discrete linear time-invariant model of structural motion that is central in damage identification under consideration is described by a second order differential equation:
\begin{equation}
\label{DynamicEquation}
M\ddot{\vec{u}}(t) + C\dot{\vec{u}}(t) + K\vec{u}(t) = \vec{f}(t),
\end{equation}
where the matrices $M$, $C$ and $K$ are real time-independent square system mass, damping and stiffness matrices of order
$d\times d$ with $d$ corresponding to the number of degrees of freedom of the model and $\vec{u}(t)$ is a time dependent displacement vector with $d$ entries.
Dots represent derivatives with respect to time $t$ and $\vec{f}(t)$ is a vector of external forces.
Considering the free vibration case, i.e. $\vec{f}(t)=\vec{0}$ and looking for the harmonic solution of Equation \eqref{DynamicEquation} in the form $\vec{u}(t)=\vec{\phi}_k\e ^{j\omega_k t}$, we obtain the following generalized eigenvalue problem
\begin{equation}
\label{EigenvalueProblem}
\left(-\omega_k^2 M +j\omega_k C + K \right)\vec{\phi}_k=\vec{0}.
\end{equation}
Here, $j=\sqrt{-1}$, $\lambda_k = \omega_k^2 = (2\pi f_k)^2$ and $\vec{\phi}_k$ are the $k^{th}$ eigenvalue and eigenvector, respectively, whereas $f_k$ is the $k^{th}$ eigenfrequency.
From Equation \eqref{EigenvalueProblem} it is easy to see that changes in system matrices $M$, $C$ and $K$ cause changes in the modal parameters $\lambda_k$ and $\vec{\phi}_k$.

It is very popular to update system matrices by the substructure
matrices \cite{Natke88SHM,FrMo95SHM,Link01SHM} as follows
\begin{equation}
\begin{aligned}
K(\bm{ \alpha}) &= K^0 - \sum_{i=1}^I\alpha_iK_i,\\
M(\bm{ \beta}) &= M^0 - \sum_{j=1}^J\beta_jM_j,\\
C(\bm{ \gamma}) &= C^0 - \sum_{s=1}^S\gamma_sC_s,
\end{aligned}
\label{SubmodelUpdating1}
\end{equation}
where $K(\vec{\alpha})$, $M(\vec{\beta})$ and $C(\vec{\gamma})$ are the improved matrices of the parameterized or corrected model. $K_i$, $M_j$ and $C_s$ are the constant expanded order matrices for the $i^{th}$, $j^{th}$ and $s^{th}$ element or
substructure (group) representing the uncertain model property and location. $\alpha_i$, $\beta_j$ and $\gamma_s$ are dimensionless
 updating parameters which can be taken as the negative relative difference of the physical parameter from its initial value, i.e. $\frac{X_t^0-X_t}{X_t^0}$, where
$t$ is one of $i$, $j$ or $s$.  This choice of updating parameters comes naturally from the simple isotropic damage theory
\cite{LemDesm05SHM}. In this theory the damage is described by a reduction in bending stiffness, as
\begin{equation}
\label{DI}
DI=\frac{E^0-E}{E^0},
\end{equation}
where $E^0$ and $E$ is the initial (undamaged) and updated (damaged) elasticity modulus, respectively, and
DI stands for damage index. The matrices $K^0$, $M^0$ and $C^0$ in \eqref{SubmodelUpdating1} are interpreted as the initial analytical system matrices
or matrices corresponding to the undamaged structure in the content of damage identification.
The model is modified only by the updating parameters for the substructure matrices.

Thus, using the simple damage model \eqref{DI} for an undamped structure whose mass does not change significantly in the degradation process,
the finite element model is parameterized by
\begin{equation}
\begin{aligned}
K(\bm{ \alpha}) &= K^0 -\sum_{i=1}^I\alpha_iK_i,\text{ where }\alpha_i=\frac{E_i^0-E_i}{E_i^0}, \\
K(\bm{ \alpha})&\phi_k(\bm{ \alpha}) =  \lambda_k(\bm{ \alpha}) M \phi_k(\bm{ \alpha}).
\end{aligned}
\label{DamageParametrization}
\end{equation}

Clearly, a small value of $\alpha_i$, or zero in the ideal case, indicates the absence of damage for a particular element or group, positive $\alpha_i$ corresponds to decrease and negative $\alpha_i$ indicates increase of the elasticity modulus for the element or group. A good damage identification method should provide positive $\alpha_i$ for the elements or groups containing damages and $\alpha_i\approx 0$ for the undamaged elements of groups.

\textbf{Remark.} The description of damage in terms of reduction in bending stiffness only is more suitable for the simple beam structures.
In the case when also torsional components of mode shapes are involved in the measurement data, it is even more advantageous to describe damage by reduction in both bending $EI$ and torsional stiffness $GI$. In the later case, one can extend the finite element model parametrization by using similar type of dimensionless parameter as for the elasticity modulus, namely
$\alpha_i^G=\frac{G_{i}^0-G_i}{G_i^0}$,
where $G_i^0$ and $G_i$ are torsional shear modulus for the initial and for the updated state, respectively. Thus, the mixed elasticity and shear modulus model parametrization will be
\begin{equation}
K(\bm{ \alpha}) = K^0 -\sum_{i=1}^I\alpha_i^E K_i^E + \alpha_i^G K_i^G, \label{DamageParametrizationBT}
\end{equation}
where $\alpha_i^E=\frac{E_i^0-E_i}{E_i^0}$, $\alpha_i^G=\frac{G_{i}^0-G_i}{G_i^0}$ and $K_i^E$ and $K_i^G$ are the nonzero parts of the element or group constant matrix $K_i$ connected to the degrees of freedom responsible for the bending and for the torsional stiffness, respectively.

\subsection{Formulation of optimization problem}
In order to solve the parameter estimation problem, we need to define so-called residual or the difference between the measured and analytical structural properties $v$, e.g. natural frequencies, mode shapes, frequency response functions (FRFs), etc. The residual is a function $\vec{r}:\mathbb{R}^n\to \mathbb{R}^m$ with $n$ corresponding to the number of updating parameters and $m$ equal to the number of
measured observations, defined by
\begin{equation}
\label{GeneralResidual}
\vec{r}(\bm \alpha) = W_\vec{v}(\vec{v}^{mea}-\vec{v}(\bm \alpha)), \footnotemark
\end{equation}
\footnotetext{Hereafter, upper index $mea$ is referring to the measured quantity. We use boldface font for vectors and $\|\cdot\|_2$ for the $l_2$ norm, i.e. $\|\vec{r}\|_2=\left( |r_1|^2+|r_2|^2+ \ldots +|r_m|^2\right)^{1/2}$. }
where $W_\vec{v}$ is a weighting matrix, which is used in order to emphasize the most significant data.

One way to minimize the difference between the measured and analytical properties is to use least squares estimation.
The objective function is then defined as the following weighted squared Euclidean norm of the residual vector:
\begin{equation}
\label{ObjectiveFunction}
f(\bm \alpha)=\frac{1}{2}(\vec{v}^{mea}-\vec{v}(\bm \alpha))^TW(\vec{v}^{mea}-\vec{v}(\bm \alpha))=\frac{1}{2}\|\vec{r}(\bm{\alpha})\|_2^2, \text{ with } W=W_\vec{v}^TW_\vec{v}.
\end{equation}
Additionally, we require that some or all updating parameters are restricted by box constraints $l_i\leq \alpha_i \leq u_i$ and thus formulate a constrained nonlinear ($\vec{r}$ depends nonlinearly on $\vec{\alpha}$) least squares problem as follows
\begin{equation}
\label{NonRegularizedProblem}
\underset{\bm{\alpha}\in \mathbb{R}^n: \bm{l}\leq \bm{\alpha} \leq \bm{u}}{\min}\frac{1}{2}\|\vec{r}(\bm \alpha)\|_2^2.
\end{equation}
The nonlinear least squares problem has no closed form solution and usually is solved by iterative methods. In the presence of noise in the measured observations, the estimated parameters found by an iterative method can have a pronounced tendency to form an oscillating pattern that makes it difficult to localize and quantify the damage (see Figures \ref{fig:Damage1-4_NoReg_13-groups} and \ref{fig:Damage1-4_NoReg_65-groups}). A standard solution of this problem is to use a regularization technique
\begin{equation}
\label{RegularizedProblem}
\underset{\bm{\alpha}\in \mathbb{R}^n: \bm{l}\leq \bm{\alpha} \leq \bm{u}}{\min}\frac{1}{2}\|\vec{r}(\bm \alpha)\|_2^2 + \lambda R(\bm \alpha),
\end{equation}
where $\lambda$ and $R$ are the regularization parameter and the regularization function, respectively. The regularization function describes the properties of the expected solution, for example, distance from the initial guess,
measure of smoothness, etc. Another type of regularization, which can be said implicitly fits the form \eqref{RegularizedProblem}
is to use the interpolation technique which was introduced in \cite{TMR02SHM}. In this paper we investigate the
regularization function $R$ being described by a total variation of the parameter vector $\bm{\alpha}$.

\subsubsection{Residuals and their derivatives}
Let us write the vector-valued residual function $\vec{r}:\mathbb{R}^n\to \mathbb{R}^m$ \eqref{GeneralResidual} in the following form:
\begin{equation}
\vec{r}(\bm \alpha) = \left( r_1(\bm \alpha),\ r_2(\bm \alpha),\ ...,\ r_m(\bm \alpha)\right)^T.
\end{equation}
Each component of $\vec{r}$ is a function $r_i: \mathbb{R}^n \to \mathbb{R}$. Moreover, the gradient $\nabla r(\bm \alpha)$, the Hessian $\nabla^2 r(\bm \alpha)$ and the Jacobian $J_\vec{r}(\bm \alpha)$ are equal to (see \cite{Datt05MB})
\begin{align}
\nabla r_i(\bm \alpha) &= \left[ \begin{array}{cccc} \frac{\partial r_i}{\partial \alpha_1} & \frac{\partial r_i}{\partial \alpha_2} & \ldots & \frac{\partial r_i}{\partial \alpha_n} \end{array}\right]^T \in \mathbb{R}^n \label{GradientRealFunction}\\
\nabla^2 r_i(\bm \alpha) &=\left[ \begin{array}{cccc}
\frac{\partial^2 r_i}{\partial^2 \alpha_1} & \frac{\partial^2 r_i}{\partial \alpha_1\partial \alpha_2} & \ldots & \frac{\partial^2 r_i}{\partial \alpha_1\partial \alpha_n} \\
\vdots & \vdots &  & \vdots\\
\frac{\partial^2 r_i}{\partial \alpha_n \partial \alpha_1} & \frac{\partial^2 r_i}{\partial \alpha_n\partial \alpha_2}& \ldots & \frac{\partial^2 r_i}{\partial^2 \alpha_n}\end{array} \right] \in \mathbb{R}^{n\times n} \label{HessianRealFunction}\\
\nabla \vec{r}(\bm \alpha) &= J_\vec{r}(\bm \alpha)^T = \left[ \begin{array}{cccc}
\frac{\partial r_1}{\partial \alpha_1} & \frac{\partial r_2}{\partial \alpha_1} & \ldots & \frac{\partial r_m}{\partial \alpha_1}\\
\vdots & \vdots &  & \vdots\\
\frac{\partial r_1}{\partial \alpha_n}& \frac{\partial r_2}{\partial \alpha_n} & \ldots & \frac{\partial r_m}{\partial \alpha_n}\end{array} \right] = \left[ \nabla r_1\ \nabla r_2\ \ldots\ \nabla r_m \right] \in \mathbb{R}^{n\times m} \label{GradientVectorValuedFunction}\\
\nabla^2 \vec{r}(\bm \alpha) &= \left[ \begin{array}{cccc}
\nabla \frac{\partial r_1}{\partial \alpha_1} & \nabla\frac{\partial r_2}{\partial \alpha_1} & \ldots & \nabla\frac{\partial r_m}{\partial \alpha_1}\\
\vdots & \vdots &  & \vdots\\
\nabla\frac{\partial r_1}{\partial \alpha_n} & \nabla\frac{\partial r_2}{\partial \alpha_n} & \ldots & \nabla\frac{\partial r_m}{\partial \alpha_n}\end{array} \right] = \left[ \nabla^2 r_1\ \nabla^2 r_2\ \ldots\ \nabla^2 r_m \right] \in \mathbb{R}^{n\times n\times m} \label{HessianVectorValuedFunction}
\end{align}

The gradient and the Hessian of $f(\bm \alpha)=\frac{1}{2}\|\vec{r}(\bm \alpha)\|_2^2=\frac{1}{2}\vec{r}(\bm \alpha)^T\vec{r}(\bm \alpha)$ are obtained by using the chain rule:
\begin{align}
\nabla f(\bm \alpha) &= \nabla \vec{r}(\bm \alpha) \vec{r}(\bm \alpha) = \sum_{j=1}^m r_j(\bm \alpha)\nabla r_j(\bm \alpha)=J_\vec{r}(\bm \alpha)^T\vec{r}(\bm \alpha) \label{Gradient_fa}\\
\nabla^2 f(\bm \alpha) &= \nabla \vec{r}(\bm \alpha) \nabla \vec{r}(\bm \alpha)^T + \nabla^2 \vec{r}(\bm \alpha) \vec{r}(\bm \alpha) \nonumber\\
&= J_\vec{r}(\bm \alpha)^TJ_\vec{r}(\bm \alpha) + \sum_{j=1}^m r_j(\bm \alpha)\nabla ^2 r_j(\bm \alpha) \approx J_\vec{r}(\bm \alpha)^TJ_\vec{r}(\bm \alpha). \label{Hessian_fa}
\end{align}
We notice here that what distinguishes the least squares from general optimization is that the second term in \eqref{Hessian_fa} for an accurate model is much less important than $J_\vec{r}(\bm \alpha)^TJ_\vec{r}(\bm \alpha)$
because the residuals are small near the solution and thus the Hessian depends only on the first-order partial derivatives of the residuals. $J_\vec{r}$ is also called the sensitivity matrix and the corresponding
finite element model updating is therefore often called as sensitivity based.

\subsubsection{Choice of residuals}
In this paper, we fit the finite element model to the data obtained by vibration tests on a reinforced concrete plate. Such experiments result in identified eigenfrequencies and mode shapes. Then, the residual is composed of two parts, the frequency residual
$\vec{r_f}(\bm \alpha)$ and the mode shape residual $\vec{r}_s(\bm \alpha)$, by $\vec{r}(\bm \alpha) = [\vec{r}_f(\bm \alpha); \vec{r}_s(\bm \alpha)]^T$.

The frequency residual $\vec{r}_f(\vec{\alpha})$ is typically a vector with entries
\footnote{The frequency residual $\vec{r}_f$ depends on the \emph{squares} of the frequencies
          $f_j(\vec{\alpha})$ and $f^{\mathrm{mea}}_{j}$.
          One possible motivation for this is that the corresponding period lengths $1/f$
          for a mass-spring system are $1/f=2\pi\sqrt{m/k}$,
          so that $(2\pi f)^2=k/m$ is a \emph{linear} function of the stiffness $k$.}
\begin{equation}
  \label{eq:FreqResiduals}
    (\vec{r}_f(\vec{\alpha}))_{j}
      \defequal
         \omega(j)\frac{\lambda^{\mathrm{mea}}_{j}-\lambda_{j}(\vec{a})}{\lambda^{\mathrm{mea}}_{j}},
    \qquad
    j=1,\ldots,m_f.
\end{equation}
where the eigenvalue $\lambda_j=\omega_j^2$ and the angular frequency $\omega_j=2\pi f_j$ correspond to the eigenfrequency $f_j$,
$m_f$ is the number of identified eigenfrequencies and $\omega(j)$ is the $j^{th}$ element of the diagonal of the weighting matrix $W_\vec{v}$.
For the undamped eigenvalue problem \eqref{EigenvalueProblem}, the eigenvalues $\lambda_j$ are all real-valued and it can be
assumed that the corresponding mode shapes are also real \cite{Bri13SHM}.
On the other hand, the measured mode shapes come from the structure with unknown damping characteristics and are usually complex.
When updating the undamped finite element model the measured complex mode shapes must be approximated with real ones \cite{FrMo95SHM}.

The
division by $\lambda^{\mathrm{mea}}_{j}$ in \eqref{eq:FreqResiduals} is done in order to obtain a similar weight for each component of
the frequency residual. Moreover, it is important
to ensure that the analytical and the measured mode shapes correspond to the same physical mode shape that is done by mode pairing, which is described below.

To define the mode shape residual, one needs to measure the similarity between two vectors. A popular choice in the literature on finite element model updating is the modal assurance criteria \cite{AlBr82SHM,All03SHM}
\begin{equation}
\label{MAC}
\text{MAC}(\vec{\phi}^{mea},\vec{\phi})=\frac{|\left<\vec{\phi}^{mea},\vec{\phi}\right>|^2}{\|\vec{\phi}^{mea}\|_2^2\|\vec{\phi}\|_2^2} = \frac{|{\vec{\phi}^{mea}}^H\vec{\phi}|^2}{\|\vec{\phi}^{mea}\|_2^2\|\vec{\phi}\|_2^2},
\end{equation}
where $\left<.,.\right>$ denotes the scalar product and $\vec{v}^H$ is the Hermitian (complex conjugate transpose) of vector $\vec{v}$. From the discussion above, it can be assumed that both the analytical and the measured mode shapes are real and
thus the Hermitian can be substituted with the transpose operator. MAC measures the difference between two vectors in terms of their collinearity and not magnitudes.
Using the MAC function one can pair analytical and measured mode shapes. For the paired mode shapes, one can then scale the mode shape $\vec{\phi}^{mea}$ to the magnitude (norm) and "orientation" of the analytical mode shape $\vec{\phi}$ by $MSF(\vec{\phi}^{mea},\vec{\phi})\vec{\phi}^{mea}$ using so-called modal
scale factor
\begin{equation}
MSF(\vec{\phi}^{mea},\vec{\phi}) = \frac{\left< \vec{\phi}^{mea},\vec{\phi}\right>}{\|\vec{\phi}^{mea}\|_2^2}.
\end{equation}
Then, the relation between norms can be checked by using the Cauchy-Schwarz inequality
\begin{align*}
\|\frac{\left< \vec{\phi}^{mea},\vec{\phi}\right>}{\|\vec{\phi}^{mea}\|_2^2}\vec{\phi}^{mea}\| \leq \frac{\|\vec{\phi}^{mea}\|_2\|\vec{\phi}\|_2}{\|\vec{\phi}^{mea}\|_2}\frac{\|\vec{\phi}^{mea}\|_2}{\|\vec{\phi}^{mea}\|_2}\leq \|\vec{\phi}\|_2.
\end{align*}
Equality holds when $\vec{\phi}^{mea}$ and $\vec{\phi}$ are collinear.
To define the mode shape residual one can, for example, use the following formula
\begin{equation}
\label{eq:ModeResiduals}
(\vec{r}_s(\vec{\alpha}))_{j,k} \defequal \omega(m_f+(j-1)d+k)\left(\frac{{\vec{\phi}_j^{mea}}^T\vec{\phi}_j(\bm \alpha)}{\|\vec{\phi}^{mea}_j\|_2^2}\phi^{mea}_{j,k} -\phi_{j,k}(\bm \alpha)\right), \ \
    j=1,\ldots,m_f,
    \
    k=1,\ldots,d
\end{equation}
to compose $\vec{r}_s=[r_{s_1}\ r_{s_2}\ \ldots r_{s_{m_f\times d}}]^T$, where the $\omega$ term is the $(m_f+(j-1)d+k)^{th}$ element of the diagonal of the weighting matrix $W_\vec{v}$ (see Equation \eqref{GeneralResidual}), for which the first $m_f$ elements are reserved to the weights of the eigenfrequencies.

 Then, the sensitivity matrix $J_\vec{r}$ is obtained using
  \begin{subequations}
  \label{R_diff}
  \begin{align}
  \frac{\partial {r_f}_j}{\partial \alpha_i} &= \omega(j)\frac{1}{\lambda_j^{mea}}\frac{\partial \lambda_j}{\partial \alpha_i} \label{Rf_diff}\\
  \frac{\partial {r_s}_{j,k}}{\partial \alpha_i} &= \omega(m_f+(j-1)d+k)\left(\frac{{\vec{\phi}_j^{mea}}^T \frac{\partial \vec{\phi}_j}{\partial \alpha_i}}{\|\vec{\phi}^{mea}_j\|_2^2}\phi^{mea}_{j,k} -\frac{\partial \phi_{j,k}}{\partial \alpha_i}\right) \label{Rs_diff} 
\end{align}
  \end{subequations}
  The derivatives of modal data with respect to the updating parameters are computed using the Fox-Kapoor formulas \cite{FoKa68SHM}.
  In the case of the finite element model parametrization \eqref{DamageParametrization} these formulas are simplified to
 \begin{subequations}
  \label{eqs:NuPhiDerivativesViaFoxKapoor}
  \begin{align}
    \label{eq:DlambdaDae}
    \pdif{\lambda_j}{\alpha_{i}}
      =& \mbox{$\vecphi_{j}$}^{\mathrm{T}}
         \pdif{K}{\alpha_{i}}
         \vecphi_{j}
      = -\mbox{$\vecphi_{j}$}^{\mathrm{T}}
         K_{i}
         \vecphi_{j}\\
    \label{eq:DphiDae}
    \pdif{\vecphi_{j}}{\alpha_{i}}
      =& \sum_{q\neq j} \frac{\mbox{$\vecphi_{i}$}^{\mathrm{T}}
                              \pdif{K}{\alpha_{i}}
                              \vecphi_{j}}%
                             {\lambda_j-\lambda_q}
                        \vecphi_{q}
      =  \sum_{q\neq j} \frac{\mbox{$\vecphi_{i}$}^{\mathrm{T}}
                              K_{i}
                              \vecphi_{j}}%
                             {\lambda_q-\lambda_j}
                        \vecphi_{q}
  \end{align}
\end{subequations}
The number of modes in \eqref{eq:DphiDae} should be big enough to contribute to well-conditioning of the sensitivity matrix $J_r$.

\subsubsection{Problem solution}
\label{subsec:Solution}
In order to solve the optimization problem \eqref{RegularizedProblem} we use the "built-in" Matlab function fmincon \cite{MatlabOptTool}, which we supply
 with the objective function value, its gradient and Hessian on each iteration step. For the nonregularized problem \eqref{NonRegularizedProblem}, the required formulas are \eqref{ObjectiveFunction}, \eqref{Gradient_fa} and
\eqref{Hessian_fa}, the residuals are computed by \eqref{eq:FreqResiduals} and \eqref{eq:ModeResiduals} and their derivatives are found
by \eqref{R_diff} and \eqref{eqs:NuPhiDerivativesViaFoxKapoor}. 
When the regularization is involved, these formulas will be modified as it is explained in the next section.

\section{Problem regularization}
We use the notion of total variation of parameters to define the regularization function in Equation \eqref{RegularizedProblem}. The total variation is then supplied with $l_2$-norm and combination of $l_1$ and $l_2$-norms. In general, the $l_1$-norm total variation regularization results in  a piecewise constant parameter estimation which keeps
sharp jumps in the parameters in the solution if they are presented and smooths out slowly varying parameters.
On the other hand
the $l_2$-norm total variation regularization not only smooths out the slowly varying parameters but also smooths
out the sharp variations \cite{ROF92SHM,Boyd04MB,Lo12SHM}.

\subsection{Total variation}
\begin{figure}[h!]
  \leavevmode
  \centering
  \includegraphics[width=0.5\linewidth]{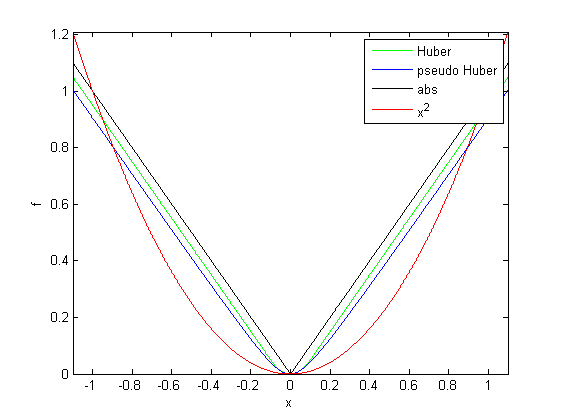}
 \caption[]
              {
                Comparison of Huber, pseudo Huber, absolute value and quadratic functions for $\mu=0.1$.
               }
      \label{fig:HuberVsPseudoHuber}
\end{figure}

Assume first that the parameter vector $\bm \alpha$ is distributed over a 2D grid as follows
\begin{equation}
\label{eq:P-p}
A=\begin{bmatrix}
 \alpha_1 &  \alpha_{d_1+1} & \cdots & \alpha_{(d_2-1)d_1+1} \\
 \alpha_2 &  \alpha_{d_1+2} & \cdots &\alpha_{(d_2-1)d_1+2} \\
  \vdots  & \vdots  & \ddots & \vdots   \\
 \alpha_{d_1} &  \alpha_{2d_1} & \cdots & \alpha_n
\end{bmatrix} \in \mathbb{R}^{d_1\times d_2}
\end{equation}
with $d_1d_2=n$, where $n$ is the number of updating parameters. Let us define the isotropic (invariant under rotations) total variation of a matrix $A$. Denote an element of this matrix at row $i$ and column $j$ by $A_{i,j}$ and define the operators
\begin{equation*}
\label{diff_operators}
\begin{aligned}[c]
D_{h_{i,j}A} = \begin{cases} A_{i+1,j}-A_{i,j}, & \mbox{if } i<d_1 \\ 0, & \mbox{if } i=d_1\end{cases}
\end{aligned}
\qquad \qquad
\begin{aligned}
D_{v_{i,j}A} = \begin{cases} A_{i,j+1}-A_{i,j}, & \mbox{if } j<d_2 \\ 0, & \mbox{if } j=d_2\end{cases}.
\end{aligned}
\end{equation*}
Compose a "discrete gradient" of $A$ by
\begin{equation*}
D_{i,jA} = \begin{bmatrix}  D_{h_{i,j}A}\\ D_{v_{i,j}A} \end{bmatrix},
\end{equation*}
where $h$ stands for differences between horizontal rows and $v$ stands for differences between vertical columns of the matrix $A$. Then, the isotropic total variation of $A$ is given by
\begin{equation}
\label{IsotropicTV}
Var_1(A) = \sum_{ij} \sqrt{D_{h_{i,j}A}^2 + D_{v_{i,j}A}^2} = \sum_{ij} \|D_{i,jA}\|_2.
\end{equation}
Now one can consider the regularized problem \eqref{RegularizedProblem} with $R(\bm \alpha)=Var_1(A)$ or $l_1$-norm total variation regularization.
Note that if $A$ is just a row or a column vector, i.e. $A=\bm \alpha$, then $Var_1(A)$ is reduced to $\sum_{j}|D_{1,jA}|$ or $\sum_{i}|D_{i,1A}|$, respectively.
Unfortunately, the function $Var_1(A)$ is not differentiable. To resolve this problem for the methods which require first order derivatives the total variation is usually modified in the following way:
\begin{equation}
\label{IsotropicTVm}
Var_\varphi(A) = \sum_{ij} \varphi \left(\|D_{i,jA}\|_2\right).
\end{equation}
The straightforward choice of $\varphi$ in Equation \eqref{IsotropicTVm} is $\varphi(x) = x^2$ and thus the $l_2$-norm total variation regularization with $R(\bm \alpha)$ in \eqref{RegularizedProblem} equals to
\begin{equation}
\label{IsotropicL2TVm}
Var_2(A)=\sum_{ij}D_{h;ijA}^2 + D_{v;ijA}^2.
\end{equation}
 However, introducing the $l_2$-norm for differentiability leads to an overregularized solution that destroys the effect of edges \cite{ROF92SHM,ChKu97SHM}
and therefore is not a good choice in cases when more precise damage localization is required. Another choice of $\varphi$ that resembles more the behavior of the absolute value function is a differentiable so-called Huber function $\varphi_\mu^H$ (\cite{Hu64Math}, Section 4, point (iii))
\begin{equation}
\label{SmooothedABS}
\varphi_\mu^{H}(x) = \begin{cases} x^2/(2\mu), & \mbox{if } |x|\leq \mu \\ |x| -\mu/2, & \mbox{if } |x| \geq \mu\end{cases}.
\end{equation}
Such defined Huber function is a smooth approximation of the absolute value function. The smaller the parameter $\mu$ the better the approximation of the absolute value function. Then, the corresponding Huber total variation is
\begin{equation}
\label{IsotropicHTV}
Var_H(A) = \sum_{ij} \varphi_\mu^H \left(\|D_{i,jA}\|_2\right).
\end{equation}
 Unfortunately, the Huber function is only first-order differentiable. Further improvement of the total variation for the second-order methods usually leads to
the so-called pseudo Huber function \cite{Har03}, which is defined by
\begin{equation}
\label{PseudoHuberfunction}
\varphi_\mu^{PH} (x) = \mu(\sqrt{1+(x/\mu)^2}-1).
\end{equation}
The pseudo Huber total variation is given by
\begin{equation}
\label{IsotropicPHTV}
Var_{PH}(A) = \sum_{ij} \phi_\mu^{PH} \left(\sqrt{D_{h_{i,j}A}^2 + D_{v_{i,j}A}^2} \right) = \sum_{ij} \varphi_\mu^{PH} \left(\|D_{i,jA}\|_2\right).
\end{equation}
For small values of $x$, the function $\varphi_\mu^{PH}$ approximates $x^2/\mu$ (use Taylor series expansion). For large values of $x$ it tends to $|x|$. It has derivatives of any order. Figure \ref{fig:HuberVsPseudoHuber} shows the difference between the Huber, pseudo Huber, absolute value and quadratic functions for $\mu=0.1$.

\subsection{$l_2$ norm total variation regularization}
Taking $\varphi(x) = x^2$ in Equation \eqref{IsotropicTVm}, we can view the problem \eqref{RegularizedProblem} as a penalized least squares problem
with regularization applied directly to problem \eqref{NonRegularizedProblem} as follows
\begin{equation}
\label{Penalized2DBV2}
\underset{\bm{\alpha}\in \mathbb{R}^n: \bm{l}\leq \bm{\alpha} \leq \bm{u}}{\min}\frac{1}{2}\|\vec{r}(\bm \alpha)\|_2^2 + \lambda Var_2(A),
\end{equation}
where $Var_2(A)$ is defined by Equation \eqref{IsotropicL2TVm} and $A$ is connected to $\bm \alpha$ by Equation \eqref{eq:P-p}.
We show now how to modify the residual vector and its Jacobian for Problem \eqref{NonRegularizedProblem}
so that the solution suggested in Section \ref{subsec:Solution} can be used for the optimization problem \eqref{Penalized2DBV2}.

Define the Toeplitz matrix $D(n)$ by
\begin{equation}
\label{ToeplitzMatrix}
 D(n) = \begin{bmatrix}
  -1 &  1 & 0 & \cdots & 0  \\
   0 & -1 & 1 & \ddots & \vdots  \\
  \vdots  & \ddots  & \ddots & \ddots & 0  \\
   0 &  \cdots & 0 & -1 & 1
 \end{bmatrix}\in \mathbb{R}^{(n-1)\times n}.
\end{equation}
Define $D_v = AD(d_2)^T$ and $D_h = D(d_1)A$. Then, $D_v$ is a $d_1\times(d_2-1)$ matrix with elements $D_{v_{i,j}A}$ and $D_h$ is a  $(d_1-1)\times d_2$ matrix with elements $D_{h_{i,j}A}$.
Taking the elements of the matrices $D_v$ and $D_h$ columnwise, we build two residual vectors $\vec{r}_v$ and $\vec{r}_h$, respectively.
It is easy to check that for $d_2>1$
\vspace{3pt}
\begin{equation*}
\vec{r}_v = \begin{bmatrix}
\alpha_{1+d_1}-\alpha_1\\
\alpha_{2+d_1}-\alpha_2\\
\cdots \\
\alpha_n-\alpha_{(d_2-1)d_1}
\end{bmatrix}
= \begin{bmatrix}
-1 &  \bovermat{$d_1-1$}{0 & \cdots & 0 &} 1 & 0 & \cdots & 0 \\
0  & -1 & 0 & \ddots  & 0 & 1 & \ddots & \vdots \\
\vdots & \ddots& \ddots & \ddots & \ddots & \ddots & \ddots & 0\\
0  &  \cdots & 0 & -1 & 0 & \cdots & 0 & 1
\end{bmatrix}
\begin{bmatrix}
\alpha_1 \\
\alpha_2 \\
\cdots \\
\alpha_n
\end{bmatrix} = J_v\bm{\alpha}
\end{equation*}
and for $d_1>1$
\begin{equation*}
\vec{r}_h=\begin{bmatrix}
\alpha_2-\alpha_1\\
\cdots \\
\alpha_{d_1}-\alpha_{d_1-1}\\
\alpha_{d_1+2}-\alpha_{d_1+1}\\
\cdots\\
\alpha_{2d_1}-\alpha_{2d_1-1}\\
\cdots\\
\alpha_{(d_2-1)d_1+2}-\alpha_{(d_2-1)d_1+1}\\
\cdots\\
\alpha_n-\alpha_{n-1}
\end{bmatrix} =
\begin{bmatrix}
   D(d_1) &     &        &   0\\
       & D(d_1) &        &   \\
       &     & \ddots &   \\
   0 & & & D(d_1)
  \end{bmatrix}
\begin{bmatrix}
\alpha_1 \\
\alpha_2 \\
\cdots \\
\alpha_n
\end{bmatrix} = J_{h}\bm{\alpha},
\end{equation*}
where $J_h$ has $d_2$ blocks $D(d_1)$ on the diagonal.
Then, the expanded residual and the Jacobian for the problem \eqref{Penalized2DBV2} are
\begin{equation*}
\vec{r} = \left[ \begin{array}{c}
\vec{r}\\
 \sqrt{2\lambda}J_v\bm{\alpha}\\
  \sqrt{2\lambda}J_h\bm{\alpha}
  \end{array} \right] \qquad \text{and} \qquad
J_\vec{r} = \left[ \begin{array}{c}
J_\vec{r}\\
\sqrt{2\lambda}J_{v} \\
\sqrt{2\lambda}J_{h}
\end{array} \right].
\end{equation*}
For the 1D problem, when for example, $d_1=1$, the matrix $D_h$ is not defined and neither are $\vec{r}_h$ or $J_{h}$. On the other hand, $J_v = D(d_2)$ and $\vec{r}_v = D(d_2)p$. Thus,
\begin{equation*}
\label{PenalizedResidual}
\vec{r} = \left[ \begin{array}{c}
\vec{r}\\
\sqrt{2\lambda}D(d_2)\bm{\alpha}
\end{array} \right] \qquad J_\vec{r} = \left[ \begin{array}{c}
J_\vec{r}\\
\sqrt{2\lambda}D(d_2)
\end{array} \right].
\end{equation*}
For Matlab code computing the updated $\vec{r}$ and $J_\vec{r}$, see the function $l2tv$ in Appendix \ref{app:TVMatlabCode}.

\subsection{Huber total variation regularization}
\label{sec:HTV}
For Equation \eqref{IsotropicTVm} with $\varphi(x) = \varphi^{H}_\mu(x)$ defined in \eqref{SmooothedABS}, we get from \eqref{RegularizedProblem} the regularization problem
\begin{equation}
\label{PenalizedHTV}
\underset{\bm{\alpha}\in \mathbb{R}^n: \bm{l}\leq \bm{\alpha} \leq \bm{u}}{\min}\frac{1}{2}\|\vec{r}(\bm \alpha)\|_2^2 + \lambda Var_{H}(A) = \underset{\bm{\alpha}\in \mathbb{R}^n: \bm{l}\leq \bm{\alpha} \leq \bm{u}}{\min} F(\bm \alpha),
\end{equation}
where $Var_{H}(A)$ is defined by Equation \eqref{IsotropicHTV} and $A$ is connected to $\bm \alpha$ by Equation \eqref{eq:P-p}. This is not a least squares problem.

In order to compute the gradient $\nabla Var_H(\vec{\alpha})$ and the Hessian $\nabla ^2 Var_H(\vec{\alpha})$, we need the first and second order derivatives of $\varphi_\mu^{H}$ with respect to the parameters $A_{i,j}$ and therefore $\alpha_k$.
For $(\varphi_\mu^{H})_{ij}=\varphi_\mu^H(\|D_{ijA}\|)$, we get
\begin{align}
\frac{\partial(\varphi_\mu^{H})_{ij}}{\partial A_{i,j}} &= \begin{cases}-\frac{D_{h_{i,j}A}+D_{v_{i,j}A}}{\mu} & \text{$\|D_{i,jA}\|_2\leq \mu$} \\
-\frac{D_{h_{i,j}A}+D_{v_{i,j}A}}{\|D_{i,jA}\|_2} &\text{$\|D_{i,jA}\|_2\geq \mu$} \end{cases} \label{eq:dh_di_j}\\
\frac{\partial(\varphi_\mu^{H})_{i-1,j}}{\partial A_{i,j}} &= \begin{cases}\frac{D_{h_{i-1,j}}}{\mu} &\text{$\|D_{i-1,jA}\|_2\leq \mu$} \\
\frac{D_{h_{i-1,j}A}}{\|D_{i-1,jA}\|_2} &\text{$\|D_{i-1,jA}\|_2\geq \mu$}
\end{cases} \label{eq:dh_di-1_j}\\
\frac{\partial(\varphi_\mu^{H})_{i,j-1}}{\partial A_{i,j}} &= \begin{cases}\frac{D_{v_{i,j-1}A}}{\mu} &\text{$\|D_{i,j-1A}\|_2\leq \mu$} \\
\frac{D_{v_{i,j-1}A}}{\|D_{i,j-1A}\|_2} &\text{$\|D_{i,j-1A}\|_2\geq \mu$} \label{eq:dh_di_j-1}
\end{cases}
\end{align}
The nonzero second-order derivatives are the following
\begin{align}
\frac{\partial^2(\varphi_\mu^{H})_{ij}}{\partial A_{i,j}^2} &= \begin{cases}\frac{2}{\mu} &\text{$\|D_{i,jA}\|_2\leq \mu$} \\
\frac{(D_{h_{i,j}A}-D_{v_{i,j}A})^2}{\|D_{i,jA}\|_2^3} &\text{$\|D_{i,jA}\|_2> \mu$} \end{cases} \label{Hd1}\\
\frac{\partial^2(\varphi_\mu^{H})_{i-1,j}}{\partial A_{i,j}^2} &= \begin{cases}\frac{1}{\mu} &\text{$\|D_{i-1,jA}\|_2\leq \mu$} \\
\frac{D_{v_{i-1,j}A}^2}{\|D_{i-1,jA}\|_2^3} &\text{$\|D_{i-1,jA}\|_2> \mu$}
\end{cases} \label{Hd2}\\
\frac{\partial^2(\varphi_\mu^{H})_{i,j-1}}{\partial A_{i,j}^2} &= \begin{cases}\frac{1}{\mu} &\text{$\|D_{i,j-1A}\|_2\leq \mu$} \\
\frac{D_{h_{i,j-1}A}^2}{\|D_{i,j-1A}\|_2^3} &\text{$\|D_{i,j-1A}\|_2> \mu$} \end{cases} \label{Hd3}\\
\frac{\partial^2(\varphi_\mu^{H})_{ij}}{\partial A_{i+1,j}\partial A_{i,j}} &=  \begin{cases} -\frac{1}{\mu} &\text{$\|D_{i,jA}\|_2\leq \mu$} \\
\frac{D_{v_{i,j}A}(D_{h_{i,j}A}-D_{v_{i,j}A})}{\|D_{i,jA}\|_2^3} &\text{$\|D_{i,jA}\|_2> \mu$} \end{cases} \label{HdKdK+1} \\
\frac{\partial^2(\varphi_\mu^{H})_{i-1,j}}{\partial A_{i-1,j}\partial A_{i,j}} &= \begin{cases}-\frac{1}{\mu} &\text{$\|D_{i-1,jA}\|_2\leq \mu$} \\
\frac{D_{v_{i-1,j}A}(D_{h_{i-1,j}A}-D_{v_{i-1,j}A})}{\|D_{i-1,jA}\|_2^3} &\text{$\|D_{i-1,jA}\|_2> \mu$}\end{cases} \label{HdKdK-1} \\
\frac{\partial^2(\varphi_\mu^{H})_{ij}}{\partial A_{i,j+1}\partial A_{i,j}} &= \begin{cases} -\frac{1}{\mu} &\text{$\|D_{i,jA}\|_2\leq \mu$} \\
\frac{D_{h_{i,j}A}(D_{v_{i,j}A}-D_{h_{i,j}A})}{\|D_{i,jA}\|_2^3} &\text{$\|D_{i,jA}\|_2> \mu$} \end{cases}  \label{HdKdK+r}\\
\frac{\partial^2(\varphi_\mu^{H})_{i,j-1}}{\partial A_{i,j-1}\partial A_{i,j}} &= \begin{cases}-\frac{1}{\mu} &\text{$\|D_{i,j-1A}\|_2\leq \mu$} \\
\frac{D_{h_{i,j-1}A}(D_{v_{i,j-1}A}-D_{h_{i,j-1}A})}{\|D_{i,j-1A}\|_2^3} &\text{$\|D_{i,j-1A}\|_2> \mu$}\end{cases} \label{HdKdK-r}\\
\frac{\partial^2(\varphi_\mu^{H})_{i-1,j}}{\partial A_{i-1,j+1}\partial A_{i,j}} &= \begin{cases}0 &\text{$\|D_{i-1,jA}\|_2\leq \mu$} \\
-\frac{D_{h_{i-1,j}A}D_{v_{i-1,j}A}}{\|D_{i-1,jA}\|_2^3} &\text{$\|D_{i-1,jA}\|_2> \mu$}\end{cases} \label{HdKdK+r-1} \\
\frac{\partial^2(\varphi_\mu^{H})_{i,j-1}}{\partial A_{i+1,j-1}\partial A_{i,j}} &= \begin{cases}0 &\text{$\|D_{i,j-1A}\|_2\leq \mu$} \\
-\frac{D_{h_{i,j-1}A}D_{v_{i,j-1}A}}{\|D_{i,j-1A}\|_2^3} &\text{$\|D_{i,j-1A}\|_2> \mu$}\end{cases} \label{HdKdK-r+1}
\end{align}
Then, the gradient $\nabla Var_H(\vec{\alpha})$ can be found by using Equations \eqref{eq:dh_di_j}--\eqref{eq:dh_di_j-1} as follows
\begin{equation*}
\nabla Var_H(\vec{\alpha})_{(j-1)d_1+i} \defequal \frac{\partial Var_{H}(\vec{\alpha})}{\partial \alpha _{(j-1)d_1+i}} =  \begin{cases}\frac{\partial(\varphi_\mu^{H})_{ij}}{\partial A_{i,j}} &\text{$i=1$, $j=1$} \\
\frac{\partial(\varphi_\mu^{H})_{ij}}{\partial A_{i,j}}+\frac{\partial(\varphi_\mu^{H})_{i-1,j}}{\partial A_{i,j}} &\text{$2\leq i \leq d_1$, $j=1$} \\
\frac{\partial(\varphi_\mu^{H})_{ij}}{\partial A_{i,j}}+\frac{\partial(\varphi_\mu^{H})_{i,j-1}}{\partial A_{i,j}} &\text{$i=1$, $2\leq j \leq d_2$}\\
\frac{\partial(\varphi_\mu^{H})_{ij}}{\partial A_{i,j}}+\frac{\partial(\varphi_\mu^{H})_{i-1,j}}{\partial A_{i,j}}+\frac{\partial(\varphi_\mu^{H})_{i,j-1}}{\partial A_{i,j}} &\text{$2\leq i \leq d_1$, $2\leq j \leq d_2$}
\end{cases}.
\end{equation*}
The Hessian $\nabla ^2 Var_H(\vec{\alpha})$ is a $(d_1d_2)\times(d_1d_2)$ symmetric matrix. So one can compute its upper triangular part $(\nabla ^2 Var_H(\vec{\alpha}))^U$ and then expand it to the symmetric matrix. Equations \eqref{Hd1}--\eqref{Hd3}
correspond to the elements on the main diagonal of the Hessian matrix. Equation \eqref{HdKdK+1} is connected to the diagonal
$[k,k+1]$, Equation \eqref{HdKdK+r} to the diagonal $[k,k+d_1]$ and Equation \eqref{HdKdK+r-1} to
 the diagonal $[k,k+d_1-1]$. Note, that Equations \eqref{HdKdK-1}, \eqref{HdKdK-r} and \eqref{HdKdK-r+1} correspond to the lower triangular part of the Hessian and thus are already considered by Equations \eqref{HdKdK+1}, \eqref{HdKdK+r} and \eqref{HdKdK+r-1}, respectively. Thus,
\begin{equation*}
\begin{aligned}
(\nabla ^2 Var_H(\vec{\alpha}))^U_{k_1,k_2} &\defequal \frac{\partial ^2 Var_{H}(\vec{\alpha})}{\partial \alpha _{k_2} \partial \alpha _{k_1}} \\
&=  \begin{cases}\frac{\partial^2(\varphi_\mu^{H})_{11}}{\partial A_{1,1}^2} &\text{\footnotesize $k_1=k_2=1$} \\
\frac{\partial^2(\varphi_\mu^{H})_{i1}}{\partial A_{i,1}^2}+\frac{\partial^2(\varphi_\mu^{H})_{i-1,1}}{\partial A_{i,1}^2} &\text{\footnotesize $k_1=k_2=i$, $2\leq i \leq d_1$} \\
\frac{\partial^2(\varphi_\mu^{H})_{1j}}{\partial A_{1,j}^2}+\frac{\partial^2(\varphi_\mu^{H})_{1,j-1}}{\partial A_{1,j}^2} & \begin{split} \text{\footnotesize $k_1$}&\text{\footnotesize $=k_2=(j-1)d_1+1$,}\\ \text{\footnotesize $2$}&\text{\footnotesize $\leq j \leq d_2$} \end{split}\\
\frac{\partial^2(\varphi_\mu^{H})_{ij}}{\partial A_{i,j}^2}+\frac{\partial^2(\varphi_\mu^{H})_{i-1,j}}{\partial A_{i,j}^2}+\frac{\partial^2(\varphi_\mu^{H})_{i,j-1}}{\partial A_{i,j}^2} & \begin{split} \text{\footnotesize $k_1$}&\text{\footnotesize $=k_2=(j-1)d_1+i$,} \\ \text{\footnotesize $2$}&\text{\footnotesize $\leq i \leq d_1$, $2\leq j \leq d_2$} \end{split}\\
\frac{\partial^2(\varphi_\mu^{H})_{ij}}{\partial A_{i+1,j}\partial A_{i,j}} & \begin{split} \text{\footnotesize $k_1$}&\text{\footnotesize $=(j-1)d_1+i$, $k_2=k_1+1$,} \\ \text{\footnotesize $1$}&\text{\footnotesize $\leq i \leq d_1-1$, $1\leq j \leq d_2$} \end{split}\\
\frac{\partial^2(\varphi_\mu^{H})_{ij}}{\partial A_{i,j+1}\partial A_{i,j}} & \begin{split} \text{\footnotesize $k_1$}&\text{\footnotesize $=(j-1)d_1+i$, $k_2=k_1+d_1$,} \\ \text{\footnotesize $1$}&\text{\footnotesize $\leq i \leq d_1$, $1\leq j\leq d_2-1$} \end{split}\\
\frac{\partial^2(\varphi_\mu^{H})_{i-1,j}}{\partial A_{i-1,j+1}\partial A_{i,j}} & \begin{split} \text{\footnotesize $k_1$}&\text{\footnotesize $=(j-1)d_1+i$, $k_2=k_1+d_1-1$,} \\ \text{\footnotesize $2$}&\text{\footnotesize $\leq i\leq d_1$, $1\leq j \leq d_2-1$} \end{split}\\
0 & \text{\footnotesize otherwise}
\end{cases}.
\end{aligned}
\end{equation*}
Then, for the objective function $F$ defined by Equation \eqref{PenalizedHTV} we have
 \begin{align*}
 F(\bm \alpha)   &= f(\bm \alpha)+\lambda Var_H(\bm \alpha),\\
 \nabla F(\bm \alpha) &= \nabla f(\bm \alpha)+\lambda \nabla Var_H(\bm \alpha),\\
 \nabla ^2 F(\bm \alpha) &= \nabla ^2 f(\bm \alpha) + \lambda \nabla ^2 Var_{H}(\bm \alpha),
 \end{align*}
 where $f(\bm \alpha)$, $\nabla f(\bm \alpha)$ and $\nabla ^2 f(\bm \alpha)$ are defined in Equations \eqref{ObjectiveFunction}, \eqref{eq:FreqResiduals}, \eqref{eq:ModeResiduals}, \eqref{Gradient_fa} and
\eqref{Hessian_fa}.

For Matlab code computing $Var_H(\bm \alpha)$, $\nabla Var_H(\bm \alpha)$ and $\nabla ^2 Var_{H}(\bm \alpha)$, see the function $htv$ in Appendix \ref{app:TVMatlabCode}.
\subsection{Pseudo Huber total variation regularization}
Equation \eqref{IsotropicTVm} with $\varphi(x)=\varphi_\mu^{PH}(x)$ defined in \eqref{PseudoHuberfunction} gives the  regularization problem
\begin{equation}
\label{PenalizedPHTV}
\underset{\bm{\alpha}\in \mathbb{R}^n: \bm{l}\leq \bm{\alpha} \leq \bm{u}}{\min}\frac{1}{2}\|\vec{r}(\bm \alpha)\|_2^2 + \lambda Var_{PH}(A) = \underset{\bm{\alpha}\in \mathbb{R}^n: \bm{l}\leq \bm{\alpha} \leq \bm{u}}{\min} F(\bm \alpha),
\end{equation}
where $Var_{PH}(A)$ is defined by Equation \eqref{IsotropicPHTV} and $A$ is connected to $\bm \alpha$ by Equation \eqref{eq:P-p}.
It is not a least squares problem. The following expressions can be used in order to find $\nabla Var_{PH}(\vec{\alpha})$ and $\nabla ^2 Var_{PH}(\vec{\alpha})$
\begin{align*}
\frac{\partial(\varphi_\mu^{PH})_{ij}}{\partial A_{i,j}} 
&= -\frac{1}{\mu} \frac{D_{h_{i,j}A}+D_{v_{i,j}A}}{\sqrt{1+\|D_{i,jA}\|_2^2/\mu^2}}\\
\frac{\partial(\varphi_\mu^{PH})_{i-1,j}}{\partial A_{i,j}}
&= \frac{1}{\mu} \frac{D_{h_{i-1,j}A}}{\sqrt{1+\|D_{i-1,jA}\|_2^2/\mu^2}}\\
\frac{\partial(\varphi_\mu^{PH})_{i,j-1}}{\partial A_{i,j}}
&= \frac{1}{\mu} \frac{D_{v_{i,j-1}A}}{\sqrt{1+\|D_{i,j-1A}\|_2^2/\mu^2}}\\
\end{align*}
and
\begin{align*}
\frac{\partial^2(\varphi_\mu^{PH})_{ij}}{\partial A_{i,j}^2} &=
\frac{1}{\mu} \frac{2+(D_{h_{i,j}A}-D_{v_{i,j}A})^2/\mu^2}{(1+\|D_{i,jA}\|_2^2/\mu^2)^{3/2}}\\
\frac{\partial^2(\varphi_\mu^{PH})_{i-1,j}}{\partial A_{i,j}^2} &= \frac{1}{\mu} \frac{1+D_{v_{i-1,j}A}^2/\mu^2}{(1+\|D_{i-1,jA}\|_2^2/\mu^2)^{3/2}} \\
\frac{\partial^2(\varphi_\mu^{PH})_{i,j-1}}{\partial A_{i,j}^2} &= \frac{1}{\mu} \frac{1+D_{h_{i,j-1}A}^2/\mu^2}{(1+\|D_{i,j-1A}\|_2^2/\mu^2)^{3/2}} \\
\frac{\partial^2(\varphi_\mu^{PH})_{ij}}{\partial A_{i+1,j}\partial A_{i,j}} &= -\frac{1}{\mu} \frac{1+D_{v_{i,j}A}(D_{v_{i,j}A}-D_{h_{i,j}A})/\mu^2}{(1+\|D_{i,jA}\|_2^2/\mu^2)^{3/2}} \\
\frac{\partial^2(\varphi_\mu^{PH})_{i-1,j}}{\partial A_{i-1,j}\partial A_{i,j}} &= -\frac{1}{\mu} \frac{1+D_{v_{i-1,j}A}(D_{v_{i-1,j}A}-D_{h_{i-1,j}A})/\mu^2}{(1+\|D_{i-1,jA}\|_2^2/\mu^2)^{3/2}} \\
\frac{\partial^2(\varphi_\mu^{PH})_{ij}}{\partial A_{i,j+1}\partial A_{i,j}} &= -\frac{1}{\mu} \frac{1+D_{h_{i,j}A}(D_{h_{i,j}A}-D_{v_{i,j}A})/\mu^2}{(1+\|D_{i,jA}\|_2^2/\mu^2)^{3/2}} \\
\frac{\partial^2(\varphi_\mu^{PH})_{i,j-1}}{\partial A_{i,j-1}\partial A_{i,j}} &= -\frac{1}{\mu} \frac{1+D_{h_{i,j-1}A}(D_{h_{i,j-1}A}-D_{v_{i,j-1}A})/\mu^2}{(1+\|D_{i,j-1A}\|_2^2/\mu^2)^{3/2}} \\
\frac{\partial^2(\varphi_\mu^{PH})_{i-1,j}}{\partial A_{i-1,j+1}\partial A_{i,j}} &= -\frac{1}{\mu} \frac{D_{h_{i-1,j}A}D_{v_{i-1,j}A}/\mu^2}{(1+\|D_{i-1,jA}\|_2^2/\mu^2)^{3/2}} \\
\frac{\partial^2(\varphi_\mu^{PH})_{i,j-1}}{\partial A_{i+1,j-1}\partial A_{i,j}} &= -\frac{1}{\mu} \frac{D_{h_{i,j-1}A}D_{v_{i,j-1}A}/\mu^2}{(1+\|D_{i,j-1A}\|_2^2/\mu^2)^{3/2}}
\end{align*}
Then, similar arguments as in Section \ref{sec:HTV} can be used to obtain $\nabla Var_{PH}(\vec{\alpha})$ and $\nabla ^2 Var_{PH}(\vec{\alpha})$ and further
 \begin{align*}
 F(\bm \alpha)   &= f(\bm \alpha)+\lambda Var_{PH}(\bm \alpha),\\
 \nabla F(\bm \alpha) &= \nabla f(\bm \alpha)+\lambda \nabla Var_{PH}(\bm \alpha),\\
 \nabla ^2 F(\bm \alpha) &= \nabla ^2 f(\bm \alpha) + \lambda \nabla ^2 Var_{PH}(\bm \alpha),
 \end{align*}
 where $f(\bm \alpha)$, $\nabla f(\bm \alpha)$ and $\nabla ^2 f(\bm \alpha)$ are defined in Equations \eqref{ObjectiveFunction}, \eqref{eq:FreqResiduals}, \eqref{eq:ModeResiduals}, \eqref{Gradient_fa} and \eqref{Hessian_fa}.

 For Matlab code computing $Var_{PH}(\bm \alpha)$, $\nabla Var_{PH}(\bm \alpha)$ and $\nabla ^2 Var_{PH}(\bm \alpha)$, see the function $phtv$ in Appendix \ref{app:TVMatlabCode}.

\subsection{Choice of $\mu$ and $\lambda$ for the total variation regularization}
The parameter $\mu$ for the Huber and the pseudo Huber function was found by testing and is approximately equal to the jump in
the elements of the parameter vector around the Damage 1 (see Section \ref{sec:Plate15Measurements}). As a rule of thumb, $\mu$ in the (pseudo) Huber controls that any variation below this value will be smoothed out
and everything above $\mu$ will be possibly kept. For this reason, $\mu$ should be chosen the same for the regularized undamaged and damaged problems.
\begin{figure}[h!]
  \leavevmode
  \centering
 \includegraphics[width=0.7\linewidth]{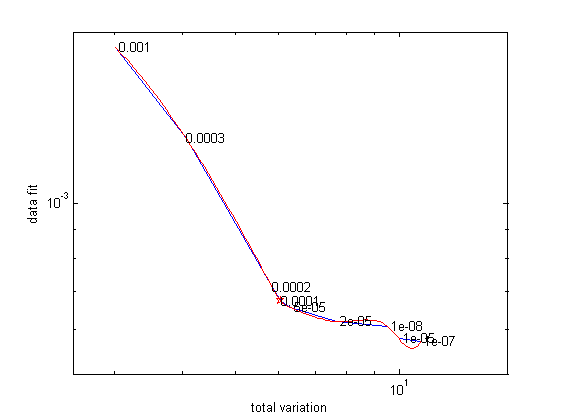}
 \caption[]
              {
               Log-log scale L-curve for the Huber total variation regularization problem with $\mu=0.01$, Damage 3 in Section
               \ref{sec:Plate15Measurements} and 65 updating parameters. The total variation corresponds to $\sum_{ij} \varphi_\mu^{H}(\|D_{ijA}\|_2)$ and data fit to $\frac{1}{2}\|\vec{r}(\bm \alpha)\|_2^2$.
               The red line is a cubic spline approximation of each straight line segment (blue). The star corresponds to the point
                on the cubic spline with maximum curvature.
               }
      \label{fig:Lcurve_damage3_htv_mu0p01}
\end{figure}
\begin{figure}[h!]
  \leavevmode
  \centering
  \includegraphics[width=0.8\linewidth]{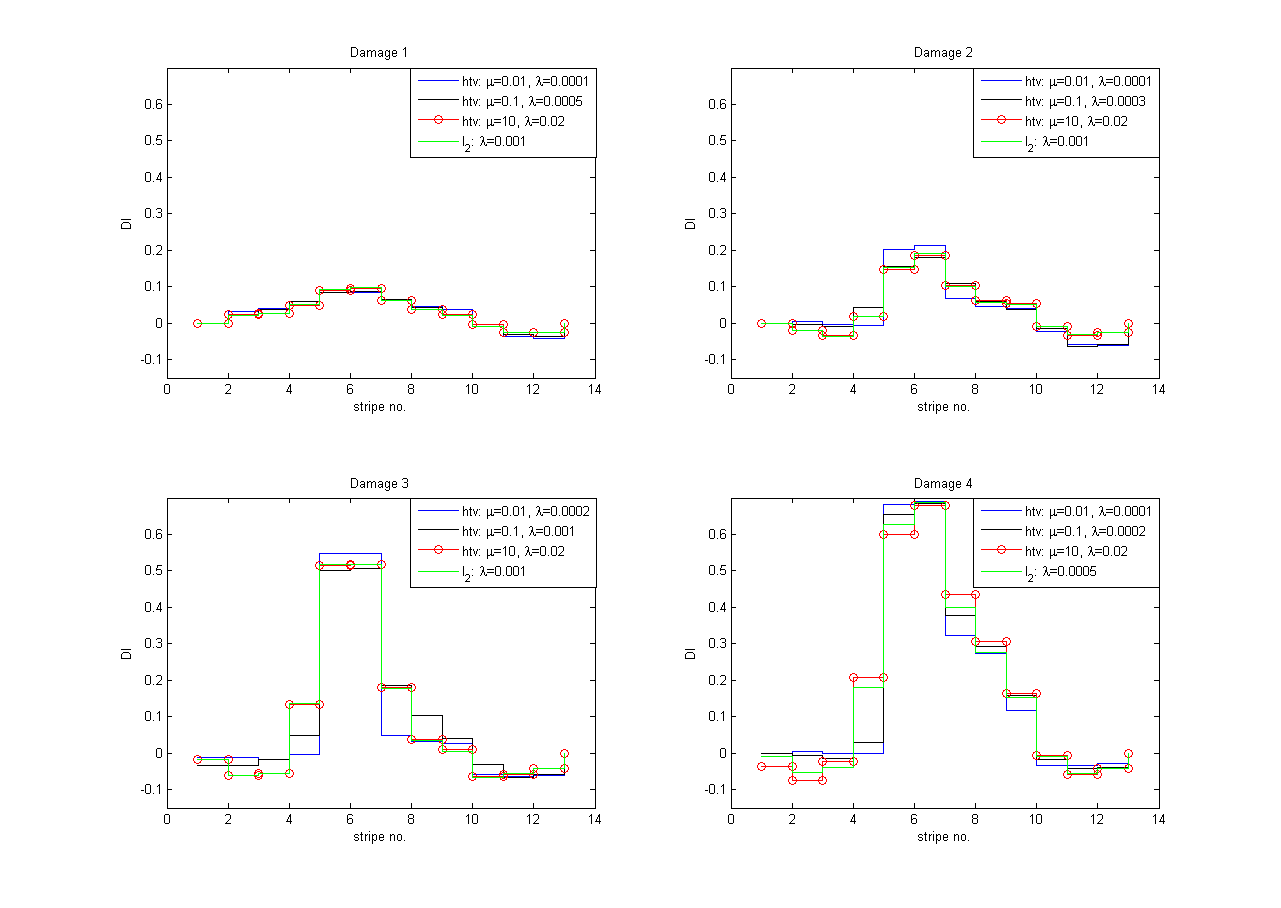}

      \caption[]
              {
                Comparison of the Huber total variation (htv) regularizations, 13 groups,
                35.9-37.2 GPa constraints for stripe no. 1 and 13 in Figure \ref{fig:GroupsOfElements}~(a), 1-40 GPa constraints for stripe no. 2-12 and
                different values of $\mu$.
               }
      \label{fig:Damage1-4_HTVDiffMu}
\end{figure}

On the other hand, to find the optimal regularization parameter $\lambda$ we use, when it is possible, the so-called L-curve method \cite{HaLe93SHM} and build a log-log-plot of the total variation norm versus
 the residual norm with $\lambda$ as a parameter. This curve shows a trade-off between doing smoothing and data fit.
In Figure \ref{fig:Lcurve_damage3_htv_mu0p01} the L-curve is drawn for Damage 3 in Section \ref{sec:Plate15Measurements} and the Huber total variation regularization with $\mu=0.01$. In the case when only a finite number of points are known on this curve,
it is popular to approximate this curve with cubic spline (red) on each line segment (blue).  Then, the
optimal $\lambda$ corresponds to the corner of the L-curve, which is defined as the point with maximal curvature of the cubic spline approximation.
 It is shown that $\lambda=0.0001$ (marked with a star) is the optimal regularization parameter. For some L-curves it was not
 possible to find the optimal value for $\lambda$ automatically. Then the parameter $\lambda$ was
chosen manually around the ``knee'' of the L-curve. A couple of different values were tested for finding one that increases
the smoothing of undesired oscillations but still keeps a sharp peak that indicates a possible damage.

More robust algorithms for finding the corner of the L-curve have been
 developed recently (e.g. \cite{CG02SHM,HaJeRo07SHM}), but here we just use a simple solution to find an estimate of the optimal $\lambda$.

Figure \ref{fig:Damage1-4_HTVDiffMu} clearly shows that when the value of $\mu$ increases, the solution with the Huber total variation regularization
becomes more similar to the solution with $l_2$-norm total variation regularization as can be expected from Equation \eqref{SmooothedABS}.

\subsection{Regularization with interpolating functions}
As it was mentioned before, another way to regularize the optimization problem \eqref{NonRegularizedProblem} is interpolation with so-called damage functions, which was suggested and used for 1D-structures in~\cite{TMR02SHM,TeDR03aSHM,TeDR05SHM,RTR10SHM}, etc. The method consists in doing
 the FEM updating with respect to the parameter vector $\vec{\alpha}$ only
for indices $p$ in a subsequence $\mathcal{P}=[\mathcal{P}_{1}, \mathcal{P}_{2}, \ldots, \mathcal{P}_{n_1}]$
of $[1,2,3,\ldots,n]$ and then use interpolation for deciding the value of the remaining parameters $\alpha_p$.

\begin{subequations}
  For example, consider a 1D-structure that is divided into 10 groups of elements with center points
  $x_p$, as illustrated in Figure~\ref{fig:2D3DtentFunctions}. The blue circles indicate a coarser grid
  of points with indices
  $\mathcal{P}=[1,4,7,10]$.
  \label{sec:InterpFunRegularization}
  \begin{figure}[h!]
    \leavevmode
    \centering
    \includegraphics[width=0.8\linewidth]{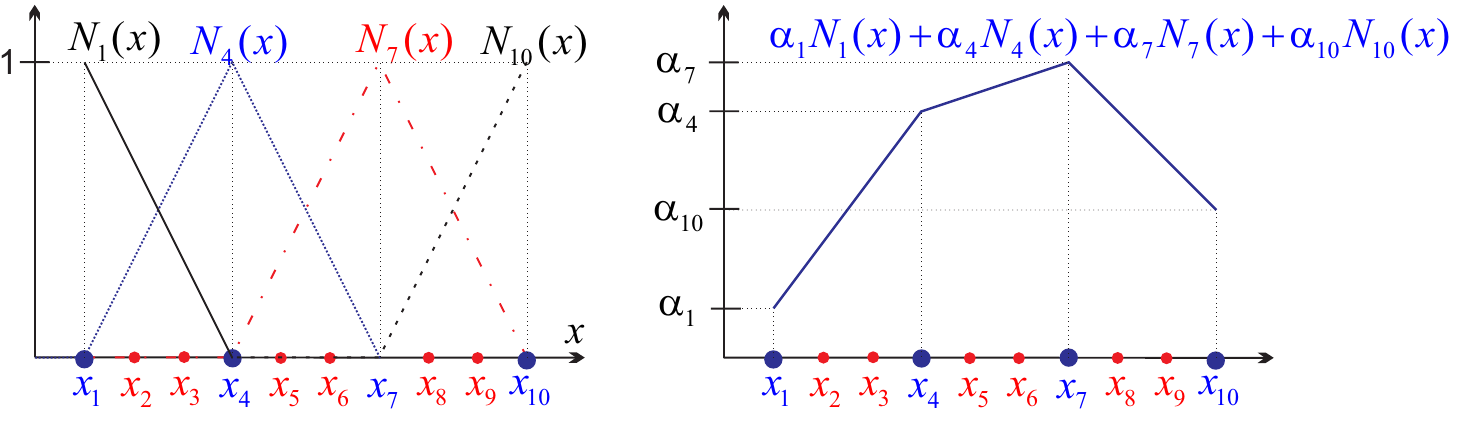}
    \caption[]
            {
              $N_1(x)-N_{10}(x)$ tent functions and their piecewise linear combination.
            }
    \label{fig:2D3DtentFunctions}
  \end{figure}

  In general, let $N_{\mathcal{P}_{k}}(x)$ be functions with the so-called
  \emph{\I{interpolation property}}
  \begin{equation}
    \label{eq:InterPolProperty}
    N_{\mathcal{P}_{k}}(x_{\mathcal{P}_{l}})
                =\delta_{l,k}
                \defequal
                 \begin{cases}
                   1 & \text{if $l=k$,}\\
                   0 & \text{if $l\neq k.$}
                 \end{cases}
  \end{equation}
  Then,
  \begin{equation}
    \label{eq:LinInterpol1D}
    \alpha(x)\defequal\sum_{l=1}^{n_1}a_{\mathcal{P}_{l}}N_{\mathcal{P}_{l}}(x)
    \qquad
    \text{and it follows that}
    \qquad
    \alpha(x_{\mathcal{P}_{l}})=\alpha_{\mathcal{P}_{l}}.
  \end{equation}
  Analogously, for a 2D-structure let $N_{\mathcal{P}_{k}}(x_{\mathcal{P}_{l}},y_{\mathcal{P}_{l}})$ be functions with the interpolation property
   \begin{equation}
    \label{eq:InterPolProperty2D}
    N_{\mathcal{P}_{k}}(x_{\mathcal{P}_{l}},y_{\mathcal{P}_{l}})
                =\delta_{l,k}
   \end{equation}
  similar to~\eqref{eq:InterPolProperty}. Again, the optimization procedure is allowed to choose the parameters $\alpha_p$ freely for
  $p\in\mathcal{P}$. Then, all other $\alpha_p = \alpha(x_p,y_p)$
  are defined as the linear interpolation by
  \begin{equation}
    \label{eq:LinInterpol2D}
    \alpha(x,y)\defequal\sum_{l=1}^{n_1}\alpha_{\mathcal{P}_{l}}N_{\mathcal{P}_{l}}(x,y) \qquad \text{ and }
    \qquad
    \alpha(x_{\mathcal{P}_{l}},y_{\mathcal{P}_{l}}) = \alpha_{\mathcal{P}_{l}} \text{ by \eqref{eq:InterPolProperty2D}}.
  \end{equation}
\end{subequations}

  The most simple choice for the 1D-case is to use piecewise linear interpolating functions
  $N_k(x)$ that are linear on each interval  $[x_{\mathcal{P}_{l}},x_{\mathcal{P}_{l+1}}]$, sometimes called
  \emph{\I{tent function}s}. For instance, see the tent functions $N_1(x)-N_{10}(x)$ and their linear combination, which indicates property \eqref{eq:LinInterpol1D}, plotted in Figure~\ref{fig:2D3DtentFunctions}.

  There are standard techniques for constructing smoother interpolating functions $N_k(x)$ with
  additional properties that are useful in signal processing~\cite[Section 3.1]{ErGr11aTF}.
  This would however not make any practical difference for the relatively small number of groups of elements used
  for the the concrete plate in Section~\ref{sec:Plate15Measurements}.

  For a 2D-structure we have tried two simple solutions. Firstly, we have generalized the interpolating functions to a 2D-grid
  similar to the suggestions in~\cite{SDY09SHM}.
  This generalization works for structures with group center points arranged in
  a rectangular grid.
  For other geometries
  further generalizations of this solution are needed. Secondly, we take the
  interpolating functions equal to \emph{\I{triangular element shape functions}} as they are defined in the FEM literature \cite{ZiTa05SHM}.
  These functions show superiority over the rectangular element interpolating functions due to much less
  restrictions on the points in the 2D grid on which they are defined. For our test case in Section~\ref{sec:Plate15Measurements}
  the results of these two solutions were pretty similar and only a bit better smoothing was achieved for the rectangular element interpolating functions used in our plots.

\begin{figure}[h!]
    \leavevmode
    \centering
    \subfloat[]{%
    \includegraphics[width=0.65\linewidth]{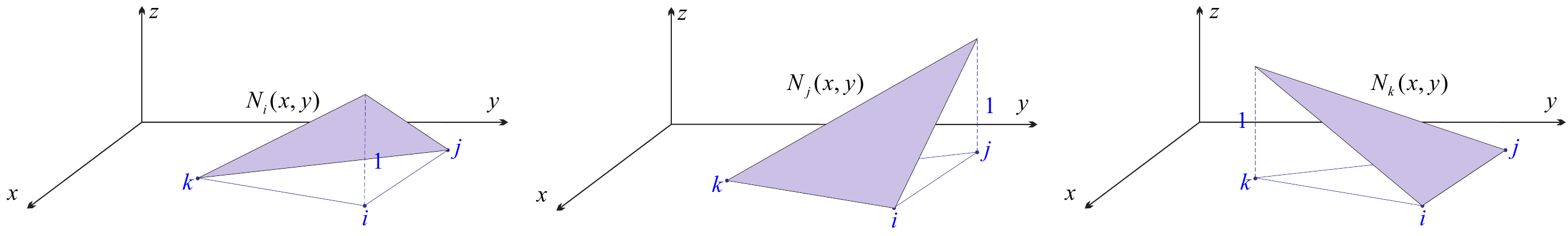}}
    \subfloat[]{%
    \includegraphics[width=0.25\linewidth]{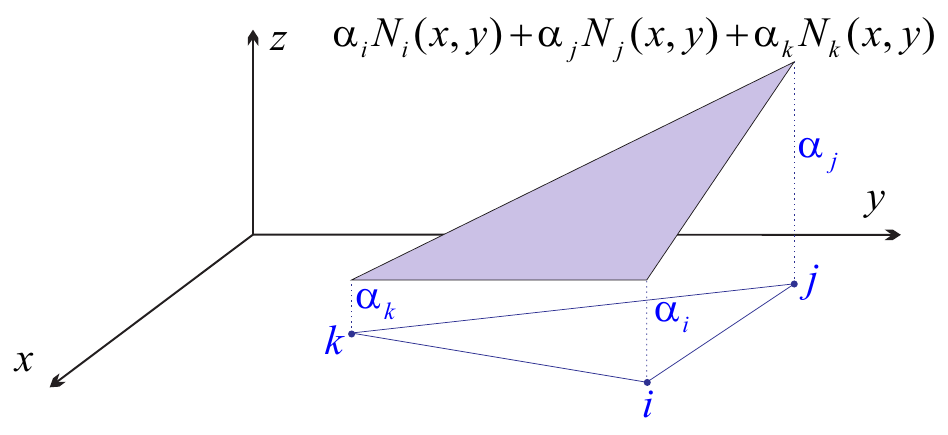}}
    \caption[]
            {
              (a) Triangular element shape functions $N_i(x,y)$, $N_j(x,y)$ and $N_k(x,y)$ for triangle $ijk$.
              (b) Their linear combination over triangle $ijk$.
            }
    \label{fig:3nodeshapeFunctions}
  \end{figure}

  Consider the case with the triangular element shape functions. Suppose that we have $P=25$ updating parameters $a_p$ in points $(x_{p},y_{p})$, ordered in a grid
  \begin{equation}
  \label{eq:2DregularGrid}
    \begin{matrix}
      \color{blue}a_{1} & \color{red}a_{ 6} & \color{blue}a_{11} & \color{red}a_{16} & \color{blue}a_{21} \\
      \color{red}a_{2} & \color{red}a_{ 7} & \color{red}a_{12} & \color{red}a_{17} & \color{red}a_{22} \\
      \color{blue}a_{3} & \color{red}a_{8} & \color{blue}a_{13} & \color{red}a_{18} &\color{blue} a_{23} \\
      \color{red}a_{4} & \color{red}a_{9} & \color{red}a_{14} & \color{red}a_{19} & \color{red}a_{24} \\
      \color{blue}a_{5} & \color{red}a_{10} & \color{blue}a_{15} & \color{red}a_{20} & \color{blue}a_{25} \\
    \end{matrix}
  \end{equation}
  with blue color for the parameters of the coarse grid.
  Now the coarse grid indices are
  \begin{equation*}
    \mathcal{P}=[1,3,5,11,13,15,21,23,25]\defequal[\mathcal{P}_{1},\mathcal{P}_{2},\ldots,\mathcal{P}_{9}].
  \end{equation*}
  We use the Delaunay triangulation \cite{MatlabFundTool}
  of the set of points $x_{\mathcal{P}_{l}}$ organized in an almost regular 2D grid similar to the one in \eqref{eq:2DregularGrid}.
  Then, for each point $(x,y)$ in the triangle $\bigtriangleup$
  with vertices at points $(x_i,y_i)$, $(x_j,y_j)$ and $(x_k,y_k)$ with area
  \begin{equation*}
  S(\bigtriangleup) = \frac{1}{2}\det \left| \begin{array}{ccc} 1 & x_i & y_i \\
                                                      1 & x_j & y_j \\
                                                      1 & x_k & y_k \end{array}\right|.
  \end{equation*}
  we define three triangular element shape functions (see Figure \ref{fig:3nodeshapeFunctions}) as follows
  \begin{align*}
  N_i^\bigtriangleup(x,y) &= \frac{1}{2S(\bigtriangleup)}\left( a_i+b_ix+c_iy\right)\\
  N_j^\bigtriangleup(x,y) &= \frac{1}{2S(\bigtriangleup)}\left( a_j+b_jx+c_jy\right)\\
  N_k^\bigtriangleup(x,y) &= \frac{1}{2S(\bigtriangleup)}\left( a_k+b_kx+c_ky\right),
  \end{align*}
  where
  \begin{equation*}
  a_i=x_jy_k-x_ky_j \qquad b_i = y_j-y_k \qquad c_i =x_k-x_j
  \end{equation*}
  and so on, with a cyclic permutation of subscripts in the order $i$, $j$ and $k$.

  \vspace{-3mm}
Let $\mat{L}$ be the $n\times n_1$-matrix
\begin{equation}
  \label{eq:interpolFctLDef}
  \mat{L}\defequal \begin{pmatrix}
                     \vec{N_{\mathcal{P}_{1}}} & \vec{N_{\mathcal{P}_{2}}} & \cdots & \vec{N_{\mathcal{P}_{n_1}}}
                   \end{pmatrix}
  \text{ with }
   \vec{N_{l}}\defequal\begin{pmatrix}
                            \tilde{N}_{l}(x_{1},y_{1})\\
                            \tilde{N}_{l}(x_{2},y_{2})\\
                            \vdots \\
                            \tilde{N}_{l}(x_{n},y_{n})
                          \end{pmatrix}
   \text{ and }
   \tilde{N}_l= \underset{\bigtriangleup_k: l\text{ vertex of }\bigtriangleup_k}{\bigcup} N_l^{\bigtriangleup_k}.
\end{equation}
If $\vec{\alpha}^{\mathcal{P}}$ is the vector of the updating parameters in  the coarse grid,
then~\eqref{eq:LinInterpol2D} and \eqref{eq:interpolFctLDef} give that
\begin{equation}
  \label{eq:aEqLcParam}
  \vec{\alpha}=\mat{L}\vec{\alpha}^{\mathcal{P}}.
\end{equation}
Moreover, the updated $m\times n_1$ Jacobian, which now can be used in Equations \eqref{Gradient_fa} and \eqref{Hessian_fa} is
\begin{equation}
\label{eq:CoarseGridJacobian}
  \mat{J}_{\vec{r}}(\vec{\alpha}^{\mathcal{P}})
            = \mat{J}_{\vec{r}}(\vec{\alpha})\mat{L},
\end{equation}
since by \eqref{eq:aEqLcParam}
\begin{align*}
  \notag
  (\mat{J}_{\vec{r}}(\vec{\alpha}^{\mathcal{P}}))_{d,l}
    \defequal&\pdif{r_d}{\alpha_{\mathcal{P}_l}}
            = \sum_{p=1}^{n}
                \pdif{r_d}{\alpha_p}
                \pdif{\alpha_p}{a_{\mathcal{P}_l}}
            = \sum_{p=1}^{n}
                (\mat{J}_{\vec{r}}(\vec{\alpha}))_{d,p}
                \pdif{(\mat{L}\vec{\alpha}^{\mathcal{P}})_p}{\alpha_{\mathcal{P}_l}}
            = \sum_{p=1}^{n}
                (\mat{J}_{\vec{r}}(\vec{\alpha}))_{d,p}
                \mat{L}_{p,l}.
\end{align*}

%
%
%
\section{Test case}
 \label{sec:Plate15Measurements}
\begin{figure}[h!]
  \leavevmode
  \centering
  \includegraphics[width=0.5\linewidth]{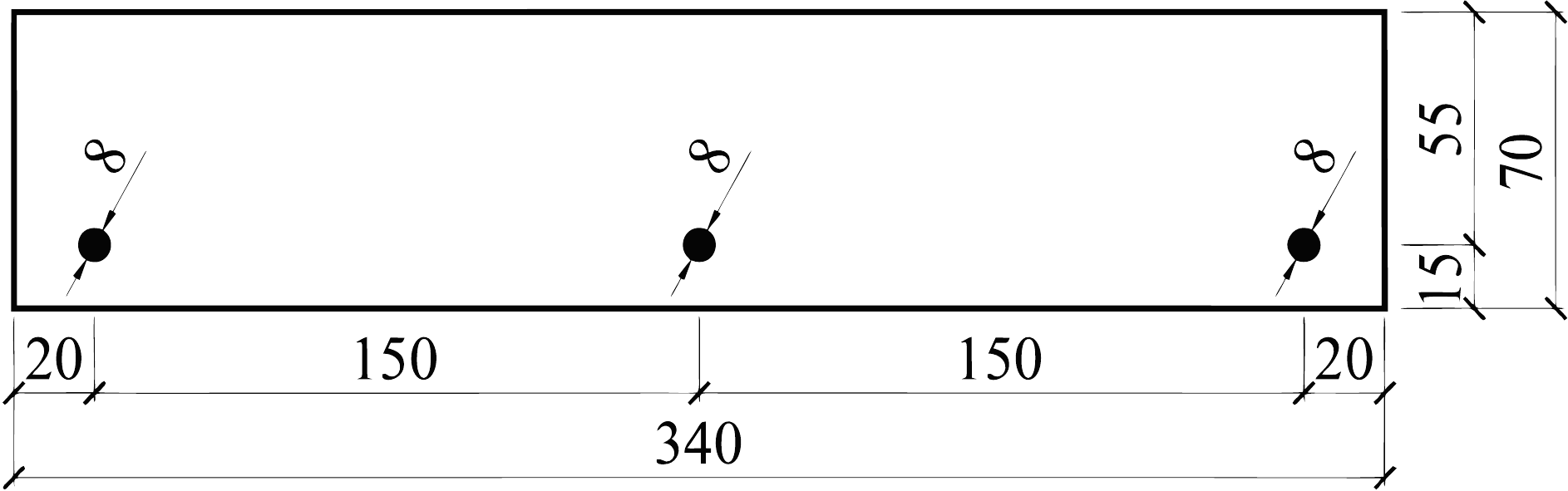}
  \vspace{3mm}
      \caption[]
              {
                Cross-section of the test plate (unit: mm).
               }
      \label{fig:TestPlateCrossSection}
\end{figure}
\begin{figure}[h!]
  \leavevmode
  \centering
    \vbox{
           \hfil
           \hbox{
                 \includegraphics[width=0.7\linewidth]{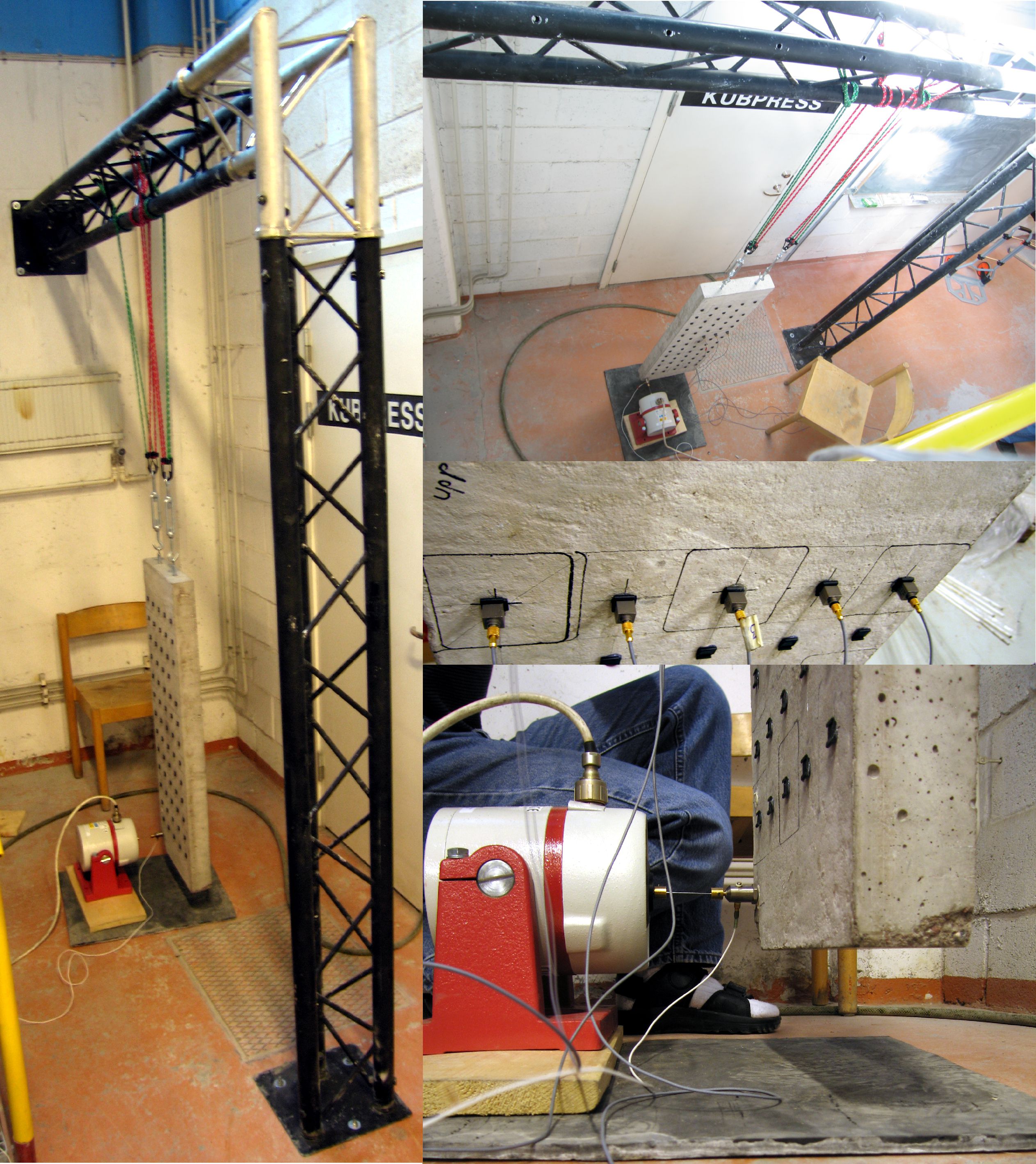}
                }
           \hfil
         }
  \caption[Measurement setup.]
          {
            The measurement setup.
          }
  \label{fig:Plate15MeasurementSetup}
\end{figure}
Measurements of forced vibrations were performed on a $1050\times 340\times 70$ mm concrete plate,
reinforced by three steel rebars of 8 mm diameter, positioned as in
Figure~\ref{fig:TestPlateCrossSection}.
The plate was excited by a swept sine force signal using an electromagnetic shaker of type LDS V406 combined with an amplifier LDS PA100E.
Ideally, either the plate or the shaker should be freely supported (or grounded) \cite[Section 3.3]{MSH97SHM}.
The plate was therefore hanging in bungee cords as shown in Figure~\ref{fig:Plate15MeasurementSetup}.
The input force and the corresponding driving point acceleration were measured by an impedance head, Br\"{u}el \& Kj{\ae}r 8001, each signal connected through a charge amplifier B\&K 2635. The remaining response points were measured using accelerometers of type B\&K 4508 B002, which were attached to the structure by using mounting plastic clips B\&K UA-1407 together with a thin layer of beeswax applied inside the clips for a more firm connection. The accelerometers/clips were glued to the plate at 5 x 13 = 65 measurement points. A B\&K 3560-C served as the data acquisition unit. It was controlled by a portable PC by using the software B\&K Pulse Labshop. These measurements were done for the following five cases
\begin{figure}[tb!]
  \leavevmode
  \centering
  \includegraphics[width=0.585\linewidth]{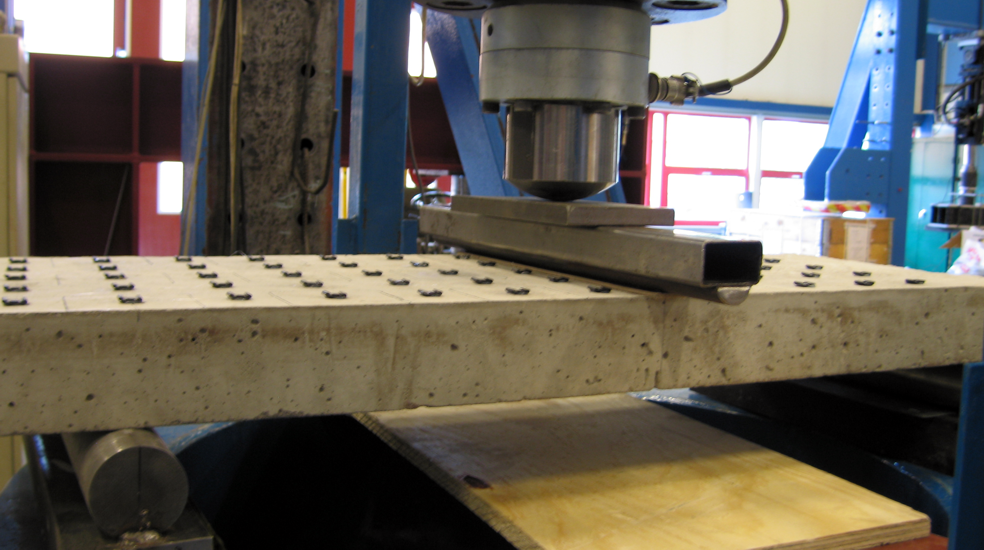}
  \includegraphics[width=0.405\linewidth]{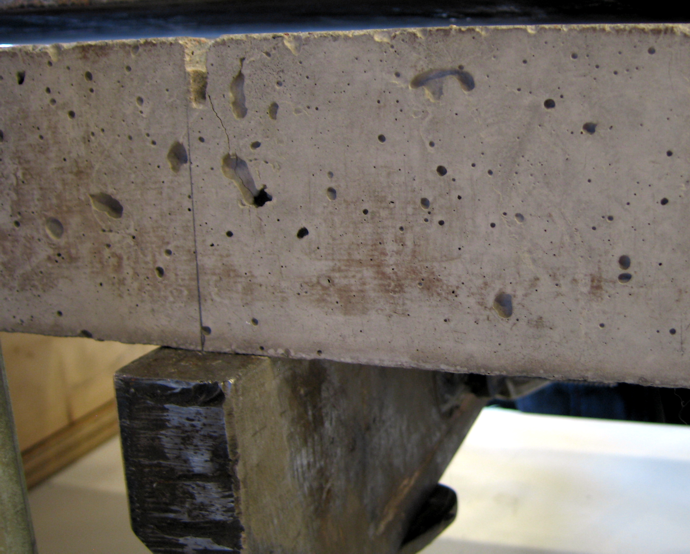}
  (a)  \hfil\hfil (b)
      \caption[]
              {
                (a) Linear 6.6 kN load for producing deeper cracks located at the notch.
                (b) Visible deeper crack after the largest applied linear load (Damage 4).
               }
      \label{fig:Cracking}
\end{figure}
\begin{description}
  \item[Damage 0] Undamaged plate.
  \item[Damage 1] A 7 mm deep notch cut with an angle grinder.
  \item[Damage 2] A 13.5 mm deep notch cut with an angle grinder.
  \item[Damage 3] Deeper real cracks, produced by applying a 6,6 kN linear load, as shown in Figure~\ref{fig:Cracking}~(a).
  \item[Damage 4] Even deeper cracks (Figure~\ref{fig:DamagePatterns}), produced by using C-clamps to apply larger linear loads.
\end{description}
\begin{figure}[h!]
  \leavevmode
  \centering
  \includegraphics[width=0.95\linewidth]{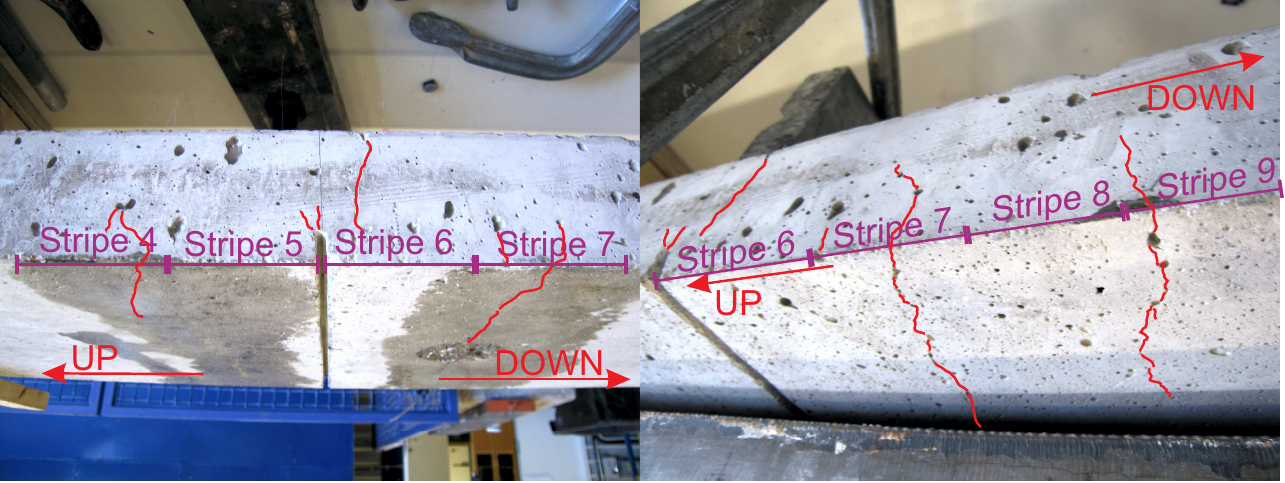}
  \\(a)  \hfil\hfil (b)
      \caption[]
              {
                Visible cracks on the plate after applying larger linear loads with C-clamps.
                The left-most crack in (a) and the right-most crack in (b) seemed to be less deep.
               }
      \label{fig:DamagePatterns}
\end{figure}
By using FRF analysis, totally 12 mode shapes were identified from the measurement data, but only the first three bending mode shapes
(see Figure \ref{fig:BendingModes135} and Table \ref{tab:eigenfreqs}) were used in the damage identification. The 30 first modes were used in the finite element analysis in order to produce system matrices
and compute  modal data derivatives with respect to the updating parameters in \eqref{eq:DphiDae}.
\begin{figure}[h!]
  \leavevmode
  \centering
  \subfloat[]{%
  \includegraphics[width=0.25\linewidth]{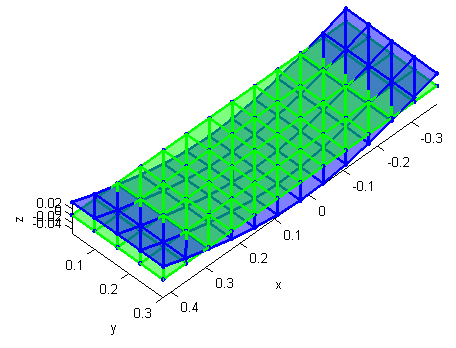}}
  \subfloat[]{%
  \includegraphics[width=0.25\linewidth]{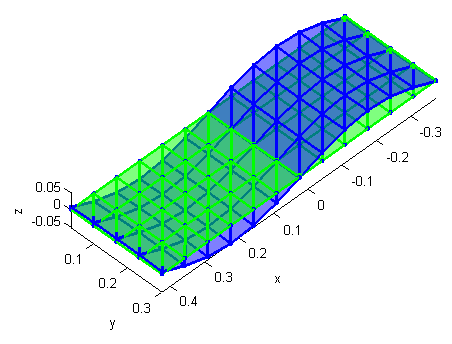}}
  \subfloat[]{%
  \includegraphics[width=0.25\linewidth]{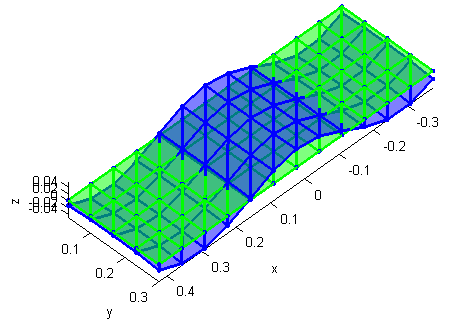}}
      \caption[]
              {
                The first three measured bending mode shapes. Green color corresponds to the undeformed plate.
                (a) Mode shape no. 1: eigenfrequency $f=249.03 \pm 0.11$ Hz.
                (b) Mode shape no. 3: $f=668.40 \pm 0.52$ Hz.
                (c) Mode shape no. 5: $f=1269.88 \pm 0.38$ Hz.
               }
      \label{fig:BendingModes135}
\end{figure}
\begin{table}[h!]
  \centering
  \begin{tabular}{|c|c|c|c|c|c|}
  \hline
    Mode & Undamaged plate & Damage 1 & Damage 2 & Damage 3 & Damage 4 \\
    \hline
    Mode 1 & $249.03 \pm 0.11$ & $243.00 \pm 0.11$ & $239 \pm 0.10$  & $217.13 \pm 0.60$  & $192.60 \pm 0.54$  \\
    Mode 3 & $668.40 \pm 0.52$  & $660.99 \pm 0.92$  & $661.48 \pm 0.22$  & $638.24 \pm 0.57$  & $604.20 \pm 0.28$  \\
    Mode 5 & $1269.88 \pm 0.38$  & $1256.97 \pm 0.80$  & $1257.13 \pm 0.94$  & $1221.70 \pm 0.58$  & $1159.06 \pm 0.53$  \\
    \hline
  \end{tabular}
  \caption{Eigenfrequencies (in Hertz with standard deviation) for the first 3 bending mode shapes for both undamaged and damaged cases.}
    \label{tab:eigenfreqs}
\end{table}

\subsection{Summary of results}
\begin{figure}[h!]
  \leavevmode
  \centering
  \subfloat[]{%
  \includegraphics[width=0.35\linewidth]{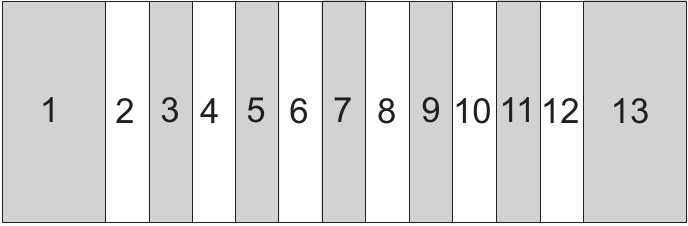}}
  \quad
  \subfloat[]{%
  \includegraphics[width=0.35\linewidth]{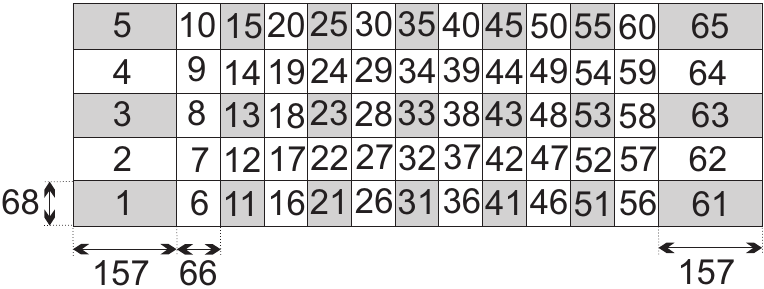}}
      \caption[]
              {
                Plate division into groups. Grey color corresponds to the coarse mesh grid used for regularization by interpolation.
                (a) 13 groups.
                (b) 65 groups (unit: mm).
               }
      \label{fig:GroupsOfElements}
\end{figure}

We have compared the results of FEMU for an 1D (wide beam) and 2D plate models and different regularization techniques.
The 1D and 2D plate models were divided into 13 and 65 groups, respectively,  as shown on Figure \ref{fig:GroupsOfElements}.

\begin{figure}[h!]
  \leavevmode
  \centering
  \includegraphics[width=0.8\linewidth]{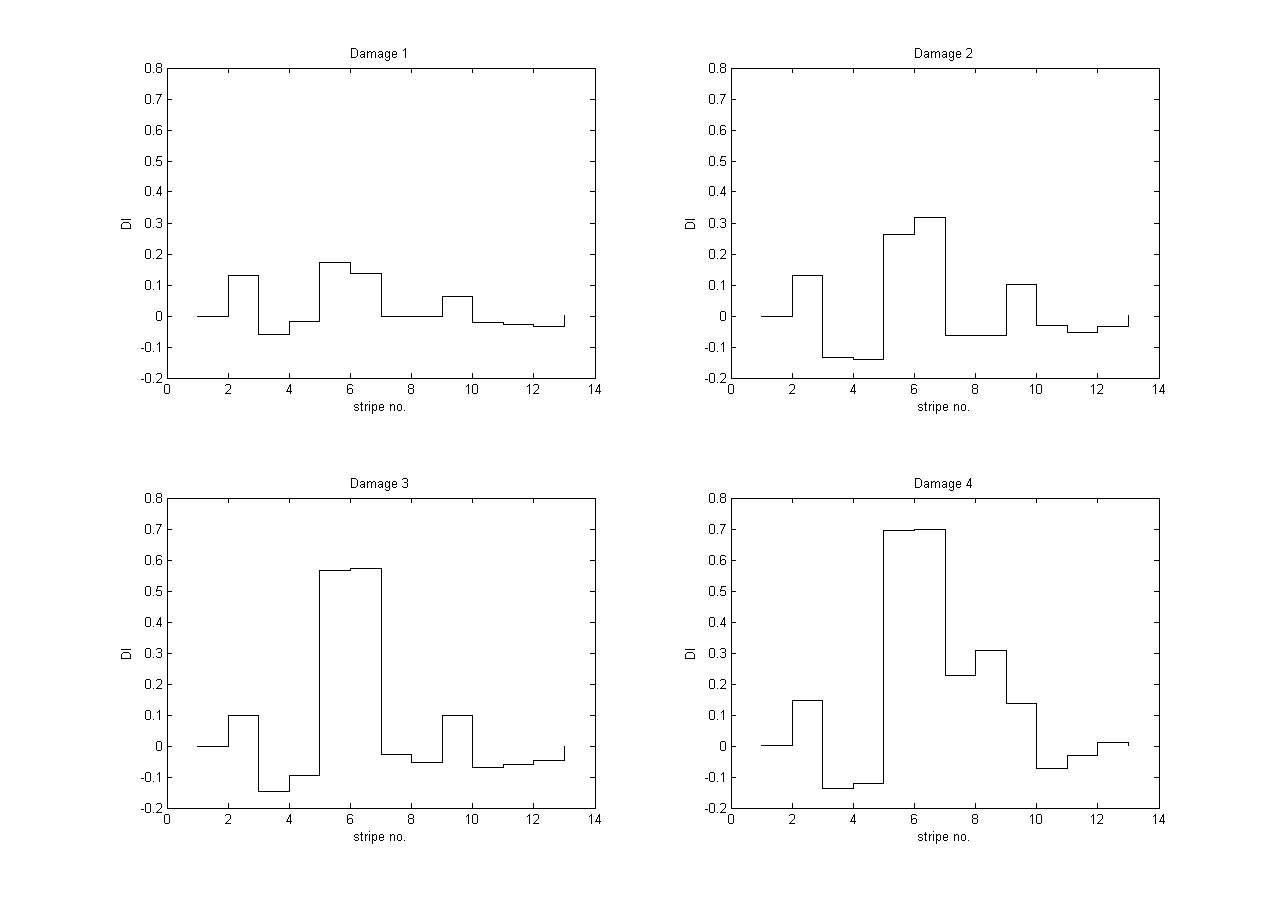}

      \caption[]
              {
                No regularization, 13 groups, 35.9--37.2 GPa constraints for the stripes no. 1 and 13, 1--40 GPa constraints
                for stripes no. 2--12.
               }
      \label{fig:Damage1-4_NoReg_13-groups}
\end{figure}
\begin{figure}[h!]
  \leavevmode
  \centering
  \includegraphics[width=0.8\linewidth]{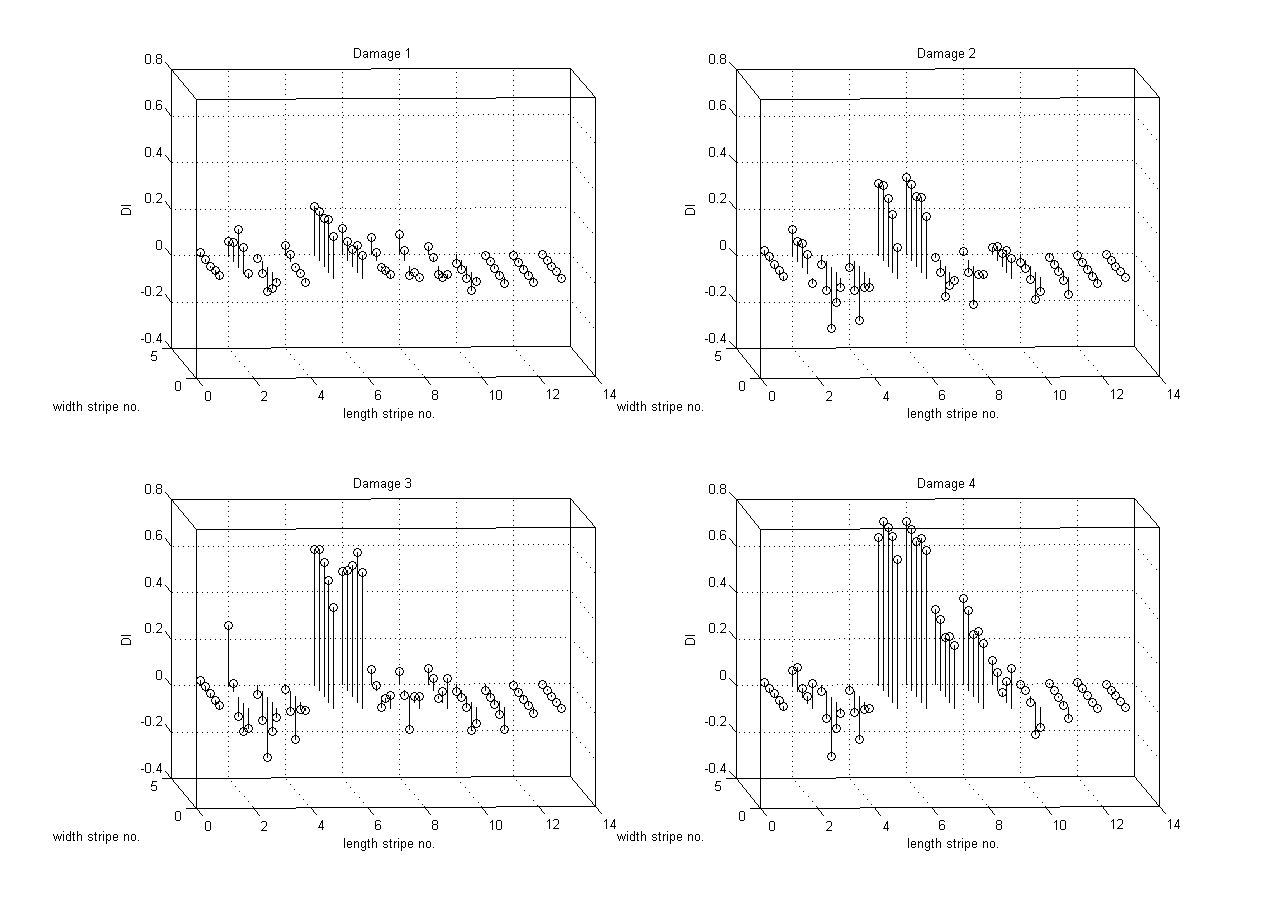}

      \caption[]
              {
                No regularization, 65 groups, 35.9--37.2 GPa constraints for short edges, 1--40 GPa constraints
                for the remaining groups.
               }
      \label{fig:Damage1-4_NoReg_65-groups}
\end{figure}
\vspace{-3mm}
The plotted parameter is the damage index, that is, $\alpha_i=DI_i=\frac{E_i^{\mathrm{0}}-E_i}{E_i^{\mathrm{0}}}$, where $E_i$ is the updated elasticity modulus and $E_i^{\mathrm{0}}$ its initial value for the $i^{th}$ group.
Thus, the  damages are indicated by high $DI_i$ values.

For all presented results, the elasticity modulus in the groups containing the two shorter edges are constrained to the range
35.9--37.2 GPa, which corresponds to maximal damage index $DI_i\approx 0.02$. This was suggested
in ~\cite{TMR02SHM} for avoiding unrealistic high
parameter values at edges,  due to lower sensitivity of modal data to changes in
elasticity modulus there.

Damages 1 and 2 (the notches) are located exactly between stripes number 5 and 6 in Figure
\ref{fig:GroupsOfElements}. 
For Damage 4, additional cracks were visible closer to and along the center line in stripe number 7, see Figure \ref{fig:DamagePatterns}.
 Thus the known damages are in the interval 5--8 and some smaller cracks in stripe 4 and 9.

The results obtained without regularization (see Figures \ref{fig:Damage1-4_NoReg_13-groups}
and \ref{fig:Damage1-4_NoReg_65-groups}) show clearly an oscillating pattern for the damage indices, from which is it
quite difficult to correctly identify both the location and severity of the damage. The damage index peaks around the
real damages, but there are also additional oscillations and peaks at stripes 2 and 9.

\begin{landscape}
\begin{figure}[h!]
  \leavevmode
  \centering
  \includegraphics[width=0.95\linewidth]{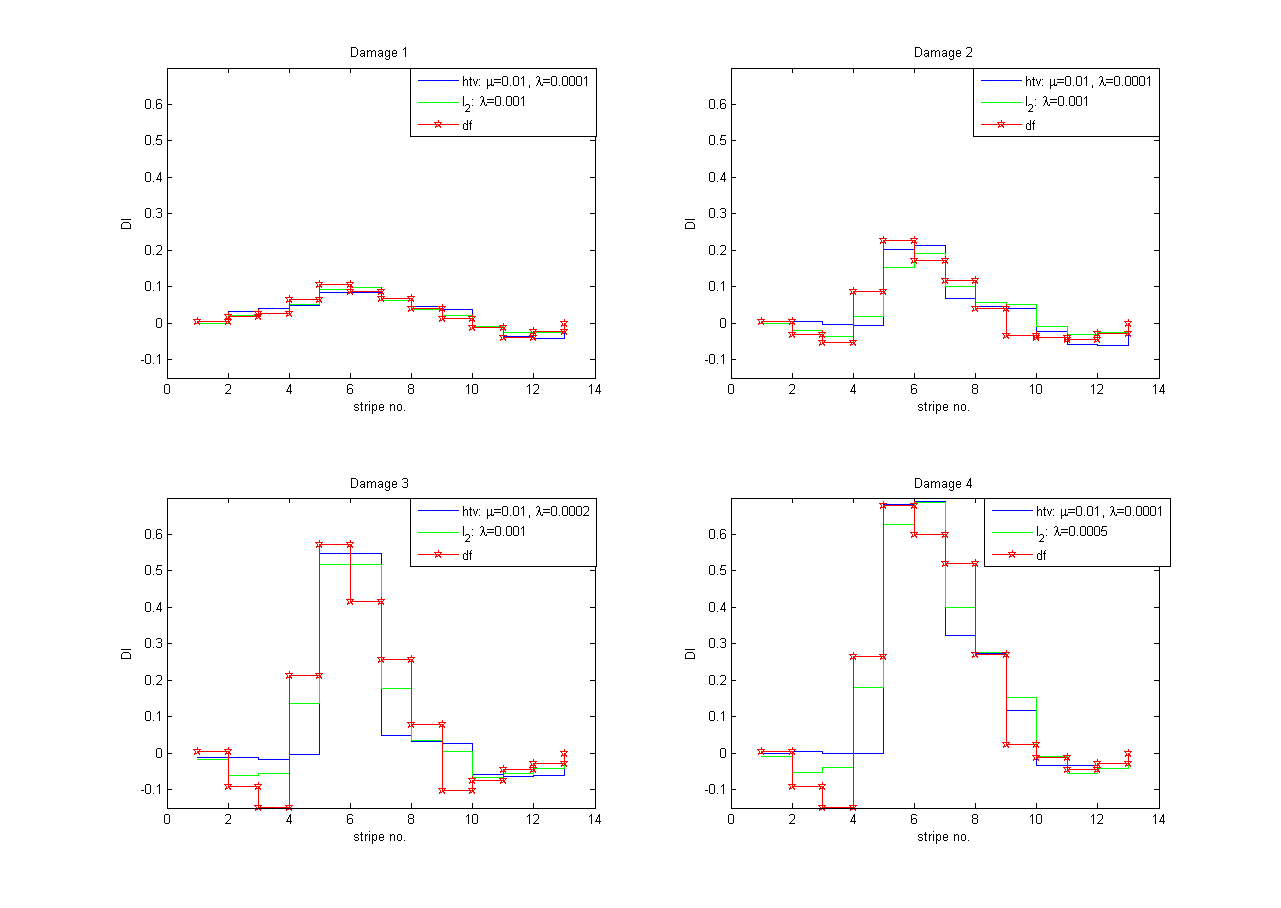}

      \caption[]
              {
                Comparison of damage functions (df), $l_2$-norm and Huber total variation (htv) regularizations, 13 groups,
                35.9-37.2 GPa constraints for stripe no. 1 and 13, 1-40 GPa constraints for stripe no. 2-12.
               }
      \label{fig:Damage1-4_13-groups_HTVvsL2vsDF}
\end{figure}
\end{landscape}

From Figure \ref{fig:Damage1-4_13-groups_HTVvsL2vsDF}, we see that for the small damage,
i.e. Damage  1 in our case, different regularization techniques result in almost the same damage pattern, which can
be described as a bell shaped parameter distribution around the damage location at position between stripes 5 and 6.
On the other hand, when the damage becomes to be more pronounced, the optimization with the Huber total variation regularization
results in a more localized damage pattern compared with the results based on either the damage functions or $l_2$-norm total
variation techniques. The fact that the test cases started with a well-localized notch (cut) damage supports this type of
damage pattern compared with more smeared bell shaped pattern.

\begin{landscape}
\begin{figure}[h!]
  \leavevmode
  \centering
  \includegraphics[width=0.95\linewidth]{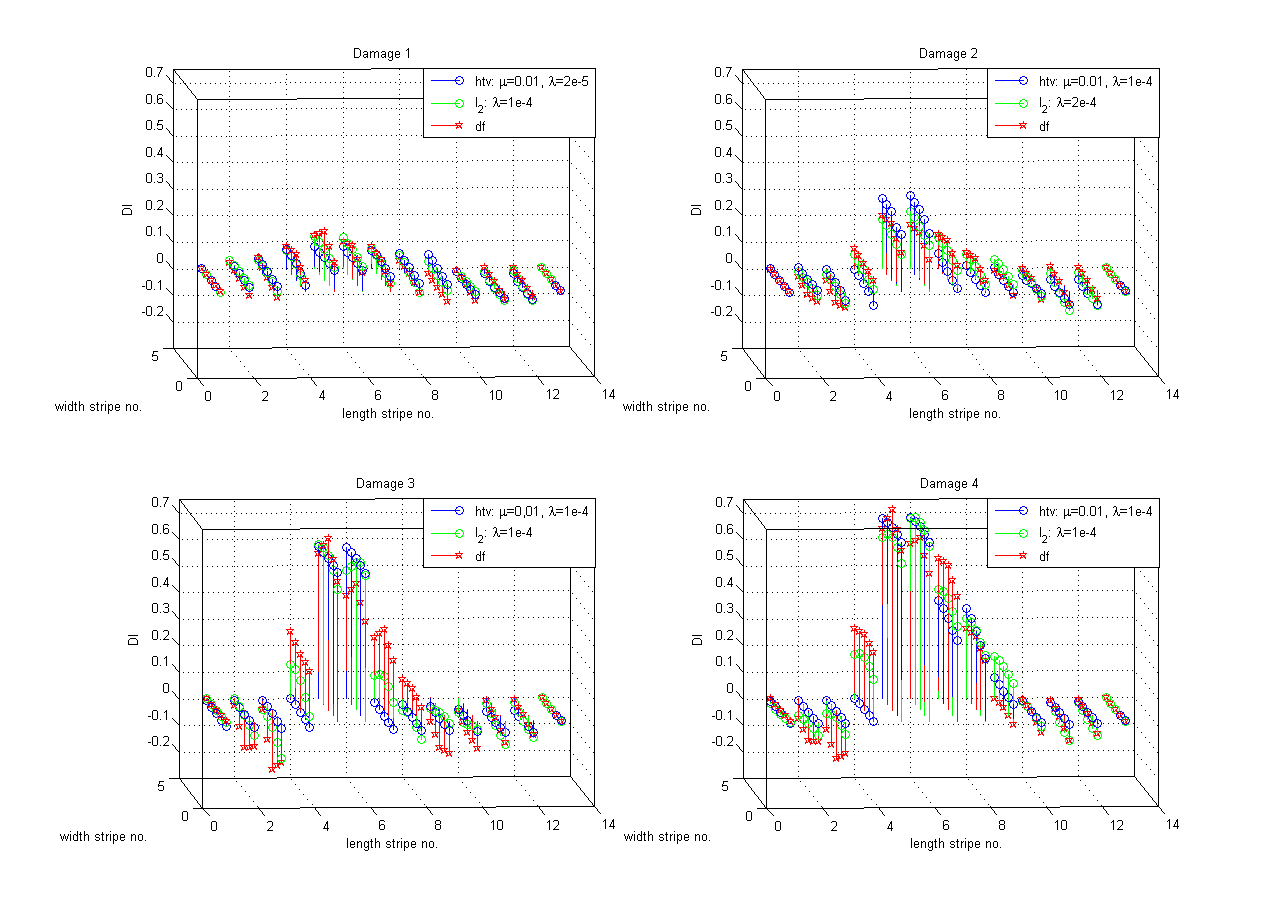}

      \caption[]
              {
                Comparison of damage functions (df), $l_2$-norm and Huber total variation (htv) regularizations, 65 groups,
                35.9-37.2 GPa constraints on short edges, 1-40 GPa constraints for the remaining groups. For the damage functions
                the coarse mesh grid is $[1\ 3\ 5\ 11\ 13\ 15\ 21\ 23\ 25\ 31\ 33\ 35\ 41\ 43\ 45\ 51\ 53\ 55\ 61\ 63\ 65]$.
               }
      \label{fig:Damage1-4_65groups_L2HTVDF}
\end{figure}
\end{landscape}

Figure \ref{fig:Damage1-4_65groups_L2HTVDF} shows almost the same comparison of methods as in Figure \ref{fig:Damage1-4_13-groups_HTVvsL2vsDF}
but for the 2D plate model and 65 groups.
Here there is a bigger difference between the results for the interpolation with the damage functions suggested in~\cite{SDY09SHM} and those for the Huber total variation.  Huber total variation gives a sharp damage indication in stripes 5--6 for Damage 2--3, corresponding
 \noindent to the notch and cracks located between these stripes, whereas for Damage 4, it also gives some indication in stripes 7--9, corresponding to the additional cracks in Figure 9 (b). It gives no indication of the small cracks in Stripe 4 in Figure 9 (b), however, which could mean that those cracks are less deep than the others. The Huber penalty term also reduces small amplitude oscillations in the updating parameters. See the differences between the methods at stripes 2--4, 5--6 and 7--12
 for Damage 2--4 in Figure~\ref{fig:Damage1-4_65groups_L2HTVDF}.

 The $l_2$-norm total variation and the damage function regularization, on the other hand,  both have the damage identification spread out over stripes 4--9 for all four damages. It is then more difficult to judge whether this indicates a spread out damage or whether it is the smoothing inherent in  these methods. For a less well-localized damage, like Damage 4, all of the above methods give more similar results.

We have also compared Huber and pseudo Huber total variation regularization. The results of these two methods are quite similar
for 65 groups and therefore are not presented here. These computations are performed in order to justify the
results obtained with the Huber total variation regularization, for which it is more easy to motivate the choice of the threshold parameter
$\mu$ but which fails to have the continuous second-order derivative required in computations.
\vspace{-3mm}
\section{Conclusions}
We have compared two different approaches for the regularization in FEMU. Interpolation based regularization,
 on the one hand, gives an
automatic smoothing of the computed updating parameters at the cost of less precise localization of the damage.
Regularization with (pseudo) Huber total variation penalty term, on the other hand, depends on a not fully automatic choice of parameters
$\mu$ and $\lambda$, but results in a more precise localization and identified severity of a well-localized damage.
An inherent advantage of the  (pseudo) Huber penalty term is that it also reduces small amplitude oscillations in the updating parameters.
 All investigated methods give more similar results for a less well-localized damage.

\vspace{-3mm}
\section{Further work}
Regularization with penalty term depends on choosing the regularization parameter $\lambda$. The L-curve and its approximation with the
cubic spline does not always give an automatic choice of the optimal $\lambda$, so a better method for finding the optimal $\lambda$ would be desired.
It would also be necessary to understand better the impact of noise on the total variation regularization methods for which
a numerical finite element model could make a contribution. It could also be interesting to apply this method to a real structure, e.g. a bridge.
\vspace{-3mm}
\section*{Acknowledgements}
We gratefully acknowledge our colleges at LTU: Inge S\"{o}derkvist for a number of valuable discussions around the regularization techniques, Fredrik Ljungren for the help with vibration tests, Lennart Elfgren, Ulf Ohlsson and Zheng Huang for the discussion about the effects of damage on the structures.

\appendix
\section{Matlab code}
\label{app:TVMatlabCode}

\begin{verbatim}
function [res,J] = l2tv(a,r,c)
% Builds a rectangular (r x c) grid A for the parameter vector a and
% then transforms it into the residual vector and the Jacobian of the
% residual at A for further use in the least squares estimation.
% This works both for 1D (r==1 or c==1) and 2D cases (r>1 and c>1).
%
% INPUT:
% a - column vector of the parameters
% r - number of rows in the 2D grid
% c - number of columns in the 2D grid
%
% OUTPUT:
% res - residual vector corresponding to A
% J - Jacobian of the residual vector at A

n = length(a);
if n~=r*c
    error('Error in grid dimension!')
end
A = reshape(a,r,c);

Dh = diff(A,[],1);
Dv = diff(A,[],2);
DDh = [];
DDv = [];

if r>1
    Dr = toeplitz([-1,zeros(1,r-2)],[-1,1,zeros(1,r-2)]);
    DDh = kron(eye(c),Dr);
end
if c>1
    DDv = toeplitz([-1,zeros(1,n-r-1)],[-1,zeros(1,r-1),1,zeros(1,n-r-1)]);
end
res = [Dv(:); Dh(:)];
J = [DDv; DDh];
\end{verbatim}
\begin{verbatim}
function [f,grad,hess] = htv(a,r,c,mu)
% Computes the value, gradient and Hessian of
% the Huber total variation at A=reshape(a,r,c).
% Huber total variation is defined as follows
% Var_phi(A) = sum_{ij} phi(sqrt((A_{i+1,j}-A_{i,j})^2 + (A(i,j+1)-A_{i,j})^2))),
% where phi is the Huber function given by
% phi(x)=x^2/(2mu) for |x|<=mu and phi(x)=|x|-mu/2 for |x|>=mu
% and mu is a predefined threshold parameter.
%
% INPUT
% a - column vector of the parameters
% r - number of rows in the 2D grid
% c - number of columns in the 2D grid, c>1
% mu - threshold parameter for the Huber function
%
% OUTPUT:
% f - value of the Huber total variation at vector a
%     organized into rectangular (r x c) grid A
% grad - gradient of Huber total variation at A
% hess - Hessian of Huber total variation at A


if c==1
  % Transpose to row vector:
  c=r; r=1;
end

A = reshape(a,r,c);

Dh = diff(A,[],1);
Dh = [Dh;zeros(1,c)];
Dv = diff(A,[],2);
Dv = [Dv zeros(r,1)];

X = sqrt(Dh.^2+Dv.^2);
F = (X<=mu).*(X.^2/(2*mu))+(X>mu).*(X-mu/2);
f = sum(sum(F));

H= (X<=mu).*Dh/mu;
V = (X<=mu).*Dv/mu;
Hs = zeros(size(H));
Vs = zeros(size(H));
Hs(X>mu) = Dh(X>mu)./X(X>mu);
Vs(X>mu) = Dv(X>mu)./X(X>mu);

Y = X.^3;
H2 = zeros(size(H));
V2 = zeros(size(H));
HV = zeros(size(H));
C = (X<=mu).*1/mu;
H2(X>mu) = Dh(X>mu).^2./Y(X>mu);
V2(X>mu) = Dv(X>mu).^2./Y(X>mu);
HV(X>mu) = Dh(X>mu).*Dv(X>mu)./Y(X>mu);

if r>1
    H0r(2:r,1:c) = H(1:r-1,1:c);
    Hs0r(2:r,1:c) = Hs(1:r-1,1:c);
    C0r(2:r,1:c)   = C(1:r-1,1:c);
    V20r(2:r,1:c)  = V2(1:r-1,1:c);
    HV0r(2:r,1:c)  = HV(1:r-1,1:c);
else
    H0r=0;Hs0r=0;C0r=0;V20r=0;HV0r=0;
end
if c>1
    V0c(1:r,2:c) = V(1:r,1:c-1);
    Vsq0c(1:r,2:c) = Vs(1:r,1:c-1);
    C0c(1:r,2:c)   = C(1:r,1:c-1);
    H20c(1:r,2:c)  = H2(1:r,1:c-1);
else
    V0c=0;Vsq0c=0;C0c=0;H20c=0;
end

grad = -H-V-Hs-Vs+H0r+Hs0r+V0c+Vsq0c;
grad = grad(:);

% (k,k) main diagonal of the Hessian
tmp = 2*C+H2-2*HV+V2+C0r+V20r+C0c+H20c;
hess = diag(tmp(:));

% (k,k+1) diagonal of the Hessian
tmp2 = -C+HV-V2;
tmp2(r,:) = zeros(1,c);
hess(r*c+1:r*c+1:end) = hess(r*c+1:r*c+1:end)+tmp2(1:r*c-1);

% (k,k+r) diagonal of the Hessian
tmp3 = -C+HV-H2;
tmp3(:,c) = []; % size(tmp3) = (r,c-1)
hess(r*c*r+1:r*c+1:end) = hess(r*c*r+1:r*c+1:end)+tmp3(1:end);

% (k,k+r-1) diagonal of the Hessian
tmp4 = -HV0r;
tmp4(:,c) = zeros(r,1);
hess(r*c*(r-1)+1:r*c+1:end) = hess(r*c*(r-1)+1:r*c+1:end)+tmp4(1:r*c-(r-1));
hess = hess + triu(hess,1)';
\end{verbatim}
\begin{verbatim}
function [f,grad,hess] = phtv(a,r,c,mu)
% Computes the value, gradient and Hessian of
% the pseudo Huber total variation at A=reshape(a,r,c).
% Pseudo Huber total variation is defined as follows
% Var_phi(A) = sum_{ij} phi(sqrt((A_{i+1,j}-A_{i,j})^2 + (A(i,j+1)-A_{i,j})^2))),
% where phi is the pseudo Huber function given by
% phi(x)=mu(sqrt(1+(x/mu)^2)-1) and mu is a predefined threshold parameter.
%
% INPUT
% a - column vector of the parameters
% r - number of rows in the 2D grid
% c - number of columns in the 2D grid, c>1
% mu - threshold parameter for the pseudo Huber function
%
% OUTPUT:
% f - value of the pseudo Huber total variation at vector a
%     organized into rectangular (r x c) grid A
% grad - gradient of pseudo Huber total variation at A
% hess - Hessian of pseudo Huber total variation at A

if c==1
  % Transpose to row vector:
  c=r; r=1;
end
A = reshape(a,r,c);

Dh = diff(A,[],1);
Dh = [Dh;zeros(1,c)];
Dv = diff(A,[],2);
Dv = [Dv zeros(r,1)];

f = sum(sum(mu*(sqrt(1+(Dh.^2+Dv.^2)/mu^2)-1)));

X  = sqrt(1+(Dh.^2+Dv.^2)/mu^2);
H  = 1/mu*Dh./X;
V  = 1/mu*Dv./X;
Y = mu*X.^3;
H2 = 1/mu^2*Dh.^2./Y;
V2 = 1/mu^2*Dv.^2./Y;
HV = 1/mu^2*Dh.*Dv./Y;
C  = 1./Y;

if r>1
    H0r(2:r,1:c)  = H(1:r-1,1:c);
    C0r(2:r,1:c)  = C(1:r-1,1:c);
    V20r(2:r,1:c) = V2(1:r-1,1:c);
    D(2:r,1:c)    = -HV(1:r-1,1:c);
else
   H0r=0;C0r=0;V20r=0;D=zeros(r,c);
end
if c>1
   V0c(1:r,2:c)  = V(1:r,1:c-1);
   C0c(1:r,2:c)  = C(1:r,1:c-1);
   H20c(1:r,2:c) = H2(1:r,1:c-1);
else
    V0c=0;C0c=0;H20c=0;
end

grad = -H-V+H0r+V0c;
grad = grad(:);

% (k,k) diagonal of the Hessian
diag1 = 2*C+H2-2*HV+V2+C0r+V20r+C0c+H20c;
hess = diag(diag1(:));

% (k,k+1) diagonal of the Hessian
diag2 = -C-V2+HV;
diag2(r,:) = zeros(1,c);
hess(r*c+1:r*c+1:end) = hess(r*c+1:r*c+1:end)+diag2(1:r*c-1);

% (k,k+r) diagonal of the Hessian
diag3 = -C-H2+HV;
diag3(:,c) = []; % size(diag3) = (r,c-1)
%diag3 = diag3'; diag3 = diag3(:);
hess(r*c*r+1:r*c+1:end) = hess(r*c*r+1:r*c+1:end)+diag3(1:r*(c-1));

% (k,k+r-1) diagonal of the Hessian
diag4 = D;
hess(r*c*(r-1)+1:r*c+1:end) = hess(r*c*(r-1)+1:r*c+1:end)+diag4(1:r*c-(r-1));

hess = hess + triu(hess,1)';
\end{verbatim} 

\bibliographystyle{unsrt}%
\bibliography{AbbrJournalNames,TimeFreq,%
                               StructuralHealthMonitoring,MathBooks}
\end{document}